\documentclass[10pt,reqno]{amsart}

\usepackage{a4wide}
\usepackage{amsmath}
\usepackage{amsfonts}
\usepackage{amssymb}
\usepackage{graphicx}
\usepackage{leqno}
\usepackage[utf8]{inputenc}
\usepackage{subfigure}
\usepackage{hyperref}

\newcommand{\deq}{\stackrel{d}{=}}
\newcommand{\mean}{\mathbb E}
\newcommand{\prob}{\mathbb P}
\newcommand{\de}{{\rm d}}

\newtheorem{Lemma}{Lemma}
\newtheorem{Proposition}{Proposition}
\newtheorem{Remark}{Remark}
\newtheorem{Theorem}{Theorem}

\newtheorem{AProposition}{Proposition}[section]

\newtheorem{Assumption}{Assumption}[section]

\numberwithin{equation}{section}

\begin{document}
\begin{center}
{\large\bf
Two parallel insurance lines  with simultaneous arrivals \\
and risks correlated with inter-arrival times}\\

\vspace{0.5cm}

\end{center}
\bigskip
\centerline{ E.S.\ Badila\footnote[1]{Supported by Project 613.001.017 of the Netherlands Organisation for Scientific Research (NWO)}, O.J.\ Boxma\footnote[2]{Supported by the IAP BESTCOM Project funded by the Belgian government} \and J.A.C.\ Resing\footnotemark[2]\let\thefootnote\relax\footnote{\emph{E-mail addresses: e.s.badila@tue.nl, o.j.boxma@tue.nl, resing@win.tue.nl.}}}

\vspace{2cm}
{\setlength{\parindent}{2pt}
\address{Department of Mathematics and Computer Science, Eindhoven University of Technology, P.O. Box 513, 5600 MB Eindhoven, The Netherlands.}

\textbf{Abstract:}
We investigate an insurance risk model that consists of two reserves which receive income at fixed
rates. Claims are being requested at random epochs from each reserve and the interclaim times
are generally distributed.
The two reserves are coupled in the sense that at a claim arrival epoch, claims are being
requested from both reserves and the amounts requested are correlated.
In addition, the claim amounts are correlated with the time
elapsed since the previous claim arrival.

We focus on the probability that this bivariate reserve process survives indefinitely.
The infinite-horizon survival problem is shown to be related to
the problem of determining the equilibrium distribution of a random walk with vector-valued increments
with 'reflecting' boundary.
This reflected random walk
is actually the waiting time process in a queueing system {\em dual} to the bivariate ruin process.

Under assumptions on the arrival process and the claim amounts,
and using Wiener-Hopf factorization with one parameter,
we
explicitly determine the Laplace-Stieltjes transform of the survival function, c.q., the
two-dimensional equilibrium waiting time distribution.

Finally, the bivariate transforms are evaluated for some examples, including for proportional reinsurance, and the bivariate
ruin functions are numerically calculated using an efficient inversion scheme.
}
\vspace{1cm}


\vspace{0.2cm}

\section{Introduction}

We study a two-dimensional ruin problem for a bivariate risk reserve process in which claims are simultaneously requested from both reserves.
The amounts of two simultaneously arriving claims may be correlated, and may also be correlated with
the time elapsed since the previous claim arrival.
Under assumptions on the arrival process and the claim sizes, we explicitly determine the Laplace-Stieltjes transform (LST) of the survival function.

Studies of {\em multidimensional} risk reserve processes are scarce in the insurance literature, although results about risk measures related to such models are highly relevant both from a theoretical and a practitioner's perspective. 
Multivariate ruin problems are relevant because they give insight into the behaviour of risk measures under various types of correlations between the insurance lines. One example is presented by multiple insurance lines within the same company which are interacting with each other as they evolve in time, via, say, coupled income rates.
Another typical example is an umbrella type of insurance model, where a claim occurrence event generates multiple types of claims which may be correlated, and each type $i$ claim is paid from its corresponding reserve $R^{(i)}$, such as car insurance together with health insurance or insurance against earthquakes. Yet another class of models is related to reinsurance problems, where a claim is shared between the insurer and one or more reinsurers.

In the existing risk and insurance literature, there are not many approaches towards analyzing such complicated
multidimensional models. A first attempt to assess multivariate risk measures can be found in the paper of Sundt~\cite{Sundt} about developing multivariate Panjer recursions  which are then used to compute the distribution of the aggregate claim process, assuming simultaneous claim events and discrete claim sizes. Other approaches are deriving integro-differential equations for the various measures of risk and then iterating these equations to find numerical approximations \cite{Chan,Gong}, or computing bounds for the different types of ruin probabilities that can occur in a setting where more than one insurance line is considered (see \cite{Cai&Li} which considers multivariate phase-type claims).
It is worth mentioning that very few papers (like \cite{APP1}, \cite{Badescu} and \cite{paper2}), analytically determine, e.g., the ruin probability for insurance models with more than one reserve.

 In an attempt to solve the integro-differential equations that arise from such models, Chan et al.~\cite{Chan} derive a so called 'boundary value problem' of a Riemann-Hilbert type for the bivariate Laplace transform of the joint survival function (see~\cite{paper2} for details about such problems arising in the context of risk and queueing theory and the book \cite{BVP} for an extended analysis of similar models in queueing). However, Chan et al. \cite{Chan} do not solve this functional equation for the Laplace transform. The law of the bivariate reserve process usually considered in the above mentioned works is that of a compound Poisson process with vector-valued jumps supported on the negative quadrant in $\mathbb R^2$, conditioned to start at some positive reserve, and linearly drifting along a direction vector that belongs to the positive quadrant. In~\cite{paper2} a similar functional equation is taken as a departure point, and it is explained how one can find transforms of ruin related performance measures via solutions of the above mentioned boundary value problems. It is also shown that the boundary value problem has an {\em explicit} solution in terms of transforms, if the claim sizes are ordered. 
  A special, important case is the setting of {\em proportional reinsurance}, which was studied in Avram et al. \cite{APP1}, \cite{APP2}. There it is assumed that there is a single arrival process, and the claims are proportionally split among two reserves. In this case, the two-dimensional exit (ruin) problem becomes a {\em one}-dimensional first-passage problem above a piece-wise linear barrier. Badescu et al.~\cite{Badescu} have extended this model by allowing a dedicated arrival stream of claims into only one of the insurance lines. They show that the transform of the time to ruin of at least one of the reserve processes can be derived using similar ideas as in \cite{APP1}.

The approach we take in this paper generalizes ideas in \cite{paper1} and \cite{paper2},
and will allow us to extend those two studies. In Section \ref{section duality} we derive a similar functional equation as in \cite{paper2} for the survival function related to a 2-dimensional reserve process, but
unlike \cite{paper2} we do not assume that the claim intervals are exponentially distributed.
Furthermore we assume the claim size vector to be correlated with the time elapsed since the previous arrival.
Such a correlation is quite natural; e.g., a claim event that generates very large claims could be subjected to additional administrative/regulatory delays.
 The type of correlation between the inter-arrival time and the vector of claim sizes is an extension to two dimensions of the dependence structure studied in \cite{paper1} for a generalized Sparre-Andersen model.
It involves making a rationality assumption regarding the trivariate LST of inter-arrival time and claim size vector (Assumption \ref{assumption rational}), which extends the case where the vector with the aforementioned components has a multivariate phase type distribution (MPH).
In addition, we also make the  assumption that the claim sizes are a.s. ordered (Assumption \ref{assumption order}).
Under these assumptions, we obtain our main result: An explicit expression for the (LST of the)
two-dimensional survival function, for a large class of vectors of interclaim times and claim amounts
of both reserves.


The paper is organized in the following way. In Section \ref{Section model assumptions} we describe the model and present the main assumptions we will be working with.
Section \ref{section duality} is dedicated to some general theory for random walks in the plane; we show a useful relation between the two-dimensional risk reserve process and a version of the random walk which has the boundary of the nonnegative quadrant in $\mathbb R^2$ as an impenetrable barrier.
This relation also makes it clear that determining the survival function
is equivalent with determining the two-dimensional waiting time distribution
in a dual two-queue two-server queueing model with simultaneous arrivals of customers at both queues.
With the help of the random walk/queueing process we derive, in Section \ref{functional}, a stochastic recursion for
the LST of the finite horizon survival function.
In Section \ref{Analysis of the stochastic recursion} we resolve the stationary version of this stochastic recursion, Formula (\ref{functional identity1}). The key tool used is a one-parameter Wiener-Hopf factorization of the bivariate kernel appearing in Equation (\ref{functional identity1}). More precisely, the Wiener-Hopf factors will depend on one parameter, which is the first argument of that bivariate kernel; see Proposition \ref{Proposition W-H 1 parameter} in Section \ref{Analysis of the stochastic recursion}.
The main result, Theorem \ref{Main result}, gives the LST of the survival function, or equivalently the stationary distribution of the waiting time/reflected random walk inside the positive quadrant (see also Remark \ref{queueing remark}).

In Section \ref{Section numerics} we explain how to calculate the transform obtained in Theorem \ref{Main result} for some examples, and we numerically calculate the ruin functions/waiting time distributions using an efficient inversion algorithm of den Iseger~\cite{denIseger}.
Finally we also point out that the numerical results suggest that the ruin functions appear to be stochastically ordered for various types of correlations between inter-arrival times and claim sizes, a positive correlation leading to smaller ruin probabilities.




\vspace{0.2cm}

 \section{Model description \label{Section model assumptions}}

 Let us begin with the general assumptions on the two reserves. The reserves start with non-negative initial capital $(u^{(1)},u^{(2)})$; as long as there are no arrivals, the reserves increase linearly with positive rates $(c^{(1)},c^{(2)})$. At the $n$th claim arrival epoch, claim sizes $B_n^{(1)}$ and $B_n^{(2)}$ are respectively requested from each reserve.
The time between the ($n-1$)th and $n$th claim arrival is denoted by $A_n$.
The sequence $\{A_n,B_n^{(1)},B_n^{(2)}\}_{n\geq 1}$, is assumed to be an i.i.d.
sequence, but within a triple, $(A_n,B^{(1)}_n,B^{(2)}_n)$ are allowed to be correlated.
We will use $A$, $B^{(1)}$, $B^{(2)}$ respectively for the generic inter-arrival time and claim sizes.
In the above-described very general set-up,
the following assumption will allow us to explicitly determine the ruin/survival probabilities by using Wiener-Hopf factorization:

 \vspace{0.2cm}
 \begin{Assumption}[\textbf{On the joint transform of $\mathbf A$, $\mathbf B^{(1)}$, $\mathbf B^{(2)}$}]\label{assumption rational} The triple transform
\begin{equation} \label{analytic assumption}H(q_0,q_1,q_2):=\mean e^{-q_0 A -q_1 B^{(1)} -q_2 B^{(2)}}\end{equation}
is a rational function in $q_i$, $i=0,1,2$, i.e., it has representation $\frac{N(q_0,q_1,q_2)}{D(q_0,q_1,q_2)}$ such that $N(q_0,q_1,q_2)$ and $D(q_0,q_1,q_2)$ are polynomials in the variables $q_i$, $i=0,1,2$.
\end{Assumption}

 $N(q_0,q_1,q_2)$ and $D(q_0,q_1,q_2)$ must satisfy some conditions, because their ratio is a transform, such as,
\[\lim\limits_{|q_0|\rightarrow \infty,\, \mathcal Re\,q_0>0} H(q_0,q_1,q_2) = \mean [e^{-q_1 B^{(1)} -q_2 B^{(2)}} 1_{\{A=0\}}].\]

 We can assume without loss of generality that $A>0$ $a.s.$ Because of the above limit, this means the degree of $N$ as a polynomial in $q_0$, $N_{q_1,q_2}(q_0)$, is strictly less than the degree of $D$ as a polynomial in $q_0$: $D_{q_1,q_2}(q_0)$, for all $q_1$ and $q_2$.

\begin{Remark}The class of rational multivariate Laplace-Stieltjes transforms contains the class of LSTs of multivariate Phase-Type distributions (MPH); see \cite{Bladt:MVME}, where the rational transform class is called multivariate matrix-exponential (MME). All of the well-known classes of multivariate Phase-Type distributions are MME. Actually, all the examples we will present are MPH distributions with a specific correlation structure which are a special case of Kulkarni’s MPH* class \cite{Kulkarni:MPH*}. There is no point in restricting ourselves to any of these subclasses. We will fully exploit the algebraic representation of rational functions in order to derive explicitly the transforms of the two-dimensional survival/ruin functions.
\end{Remark}

\vspace{0.2cm}
 The reserve process $\bar R_t= (\bar R_t^{(1)},\bar R_t^{(2)})$ evolves as

 \[\bar R_t =  \bar u +  c t- \sum_{i=1}^{n(t)} B_i,\,\mbox{where } \bar u:=(\bar u^{(1)},\bar u^{(2)}),\, c:=(c^{(1)},c^{(2)}),\, B_i:=(B^{(1)}_i,B^{(2)}_i), \]
and $n(t)$ is the number of arrivals before $t$.
We use the notation $\bar{R}$ with a bar because  later on we want to scale the reserve process by dividing
$\bar{R}_t^{(i)}$ by $c^{(i)}$, using the notation $R_t^{(i)}$ without a bar for the resulting process.
In this model,
the two reserves are correlated due to simultaneous claim arrivals {\em and} due to correlations that may exist in the claim size vector $(B^{(1)},B^{(2)})$.

This paper is concerned with measuring the event that both reserve processes survive indefinitely, i.e.,
we aim to determine the {\em survival function}

\[\bar F^s(\bar u^{(1)},\bar u^{(2)}) :=\prob (\bar R^{(1)}_t\geq 0,\;\forall t>0\mbox{ AND }\bar R^{(2)}_t\geq 0,\;\forall t> 0\;|\;\bar R_0=(\bar u^{(1)},\bar u^{(2)})).\]
In terms of times to ruin, $\bar F^s$ is related to the first time at least one of the two insurance lines is ruined,

\begin{equation}\label{naive time to ruin}\bar\tau_\wedge(\bar u^{(1)},\bar u^{(2)}) =\inf \{t;\, \min (\bar R^{(1)}_t,\bar R^{(2)}_t) <0\}= \bar\tau^{(1)}(\bar u^{(1)})\wedge \bar\tau^{(2)}(\bar u^{(2)}), \end{equation}
where $\bar\tau^{(i)}(u^{(i)})$ are the marginal times to ruin, $i=1,2$.
In particular,

\[\bar{F}^s(\bar u^{(1)},\bar u^{(2)}) =1- \prob(\bar{\tau}_\wedge(\bar u^{(1)},\bar u^{(2)}) < \infty).\]
We can also define the first time at which both insurance lines are ruined:

\[\bar\tau_{\vee}(\bar u^{(1)},\bar u^{(2)}) = \bar\tau^{(1)}(\bar u^{(1)})\vee \bar\tau^{(2)}(\bar u^{(2)}).\]
It is similarly related to the probability that \emph{at least} one of the two reserves survives indefinitely,

\[ \bar F^s_{OR} (\bar u^{(1)},\bar u^{(2)}):= \prob( \bar R^{(1)}_t\geq 0,\;\forall t>0\;{\rm OR} \;\bar R^{(2)}_t\geq 0,\;\forall t> 0\;|\;\bar R_0=(\bar u^{(1)},\bar u^{(2)})),\]
by

\[\bar{F}^s_{OR}(\bar u^{(1)},\bar u^{(2)}) = 1- \prob(\bar{\tau}_\vee(\bar u^{(1)},\bar u^{(2)}) < \infty).\]
Notice that $\bar\tau_\vee$ is not the same as $\inf \{t;\, \max (\bar R^{(1)}_t,\bar R^{(2)}_t) <0\},$ that is, joint ruin may not happen simultaneously.

\vspace{.2cm}

\noindent The above survival functions are then related by

 \[   \bar F^s_{OR} (\bar u^{(1)},\bar u^{(2)}) = \bar F^s(\bar u^{(1)},\infty) + \bar F^s(\infty,\bar u^{(2)}) - \bar F^s(\bar u^{(1)},\bar u^{(2)}),\]
where $\bar F^s(\bar u^{(1)},\infty)$ and $\bar F^s(\infty,\bar u^{(2)})$ are the marginal survival functions.
Moreover, we also have  $\bar F^{i,j}(\bar u^{(1)},\bar u^{(2)})$, the probability that reserve $i$ survives indefinitely, while reserve $j$ ruins, for $i,j=1,2$, $i\neq j$.
\[\bar F^{1,2} (\bar u^{(1)},\bar u^{(2)}) = \bar F^s (\bar u^{(1)},\infty) - \bar F^s(\bar u^{(1)},\bar u^{(2)}),\]
and similarly for $\bar F^{2,1}(\bar u^{(1)},\bar u^{(2)}).$
In view of the above, it suffices to determine $\bar F^s (\bar u^{(1)},\bar u^{(2)})$ in order to obtain all
the other survival/ruin functions.

\begin{Remark}\label{Remark scaling} For the study of survival functions, we can normalize the reserve processes by their respective income rates. The survival function is preserved, with the starting capital scaled accordingly. To be more precise, let

\[ R^{(i)}_t =  \bar u^{(i)}/c^{(i)} +   t- \sum_{k=1}^{n(t)} B^{(i)}_k/c^{(i)},\,\,i=1,2.\]
Since $ R^{(i)}_t\geq  0 \Leftrightarrow \bar R_t^{(i)}\geq 0$,  for $\bar  \tau_\wedge$, $\tau_\wedge$ the exit times of $\bar R_t$ respectively $R_t$ from the non-negative quadrant, we have the relation
\[\tau_\wedge(u^{(1)},u^{(2)}) = \bar  \tau_\wedge(c^{(1)}u^{(1)},c^{(2)}u^{(2)})\]
and then also $ F^s(u^{(1)},u^{(2)}) = \bar F^s(c^{(1)}u^{(1)},c^{(2)}u^{(2)})$, with $F^s$ the survival function of the scaled process.

This means that for our purposes it suffices to study the process $R_t$ and the associated survival functions $F^s(u^{(1)},u^{(2)})$. $F^s_{OR}$ and the complementary ruin time $\tau_{\vee}$ are similarly defined.
\end{Remark}

\vspace{0.2cm}
The main idea of our approach is to exploit the fact that the embedded process at arrival epochs of claims is a random walk in the plane with increments $(A_n-B^{(1)}_n/c^{(1)},A_n-B^{(2)}_n/c^{(2)})$, conditioned on starting at $(u^{(1)},u^{(2)})$.
The different ways in which the reserve process can be ruined correspond to the possible positions of the random walk
\begin{equation}
R_{t_n}:= u+\sum_{i=1}^n X_i,\,\,\,\,\,\,\,\,X_n:= (X_n^{(1)},X_n^{(2)}) := (A_n-B_n^{(1)}/c^{(1)}, A_n-B^{(2)}_n/c^{(2)}),
\label{RtnXn}
\end{equation}
at the time of exit from the non-negative quadrant. Here and below, $t_n$ denotes the $n^{th}$ arrival epoch.

  This is a difficult model to analyze in full generality; in particular, it is much more general than the 2-dimensional ruin process described in~\cite{Chan} and \cite{paper2}.
However, in the 
case that the claims in insurance line $1$ are
larger than those in insurance line $2$, we are able to determine the two-dimensional survival function.

 \vspace{0.2cm}
 \begin{Assumption}[\textbf{Ordering assumption}]\label{assumption order} For a generic claim event $(B^{(1)},B^{(2)})$:

\begin{equation}\label{ordering assumption} B^{(1)}/c^{(1)} \geq B^{(2)}/c^{(2)}\,\,\, a.s.\end{equation}

\end{Assumption}
\vspace{0.2cm}


 A special, important example for this ordering assumption is the case when there is a single arrival process such that the common claim is partitioned into fixed proportions $(\alpha, 1-\alpha)$, and we can always take w.l.o.g. $\alpha\in [1/2,\,1]$ ($\alpha$ may even be a random variable with this interval as support). This special case then significantly generalizes the setting in~\cite{APP1}, where it is assumed that the common arrival process is compound Poisson (and in particular the inter-arrival times are independent of the claim sizes).
%

One can go a step further and assume there is a dedicated renewal-arrival stream of claims into the line which pays the greater share, say $\alpha$. This is in line with the assumption in Badescu et al. \cite{Badescu} of a dedicated Poissonian stream of claims, and extends it, once we combine it with Assumption \ref{assumption rational}. Clearly this is not a proportional reinsurance problem anymore.
From a mathematical perspective,  the analysis is more insightful than if one would just assume proportionality.

\vspace{0.2cm}
We take the following approach: using a duality argument (Section \ref{section duality}), we derive a recursive equation for the survival function (Section \ref{functional}).
Using complex function theory, and under Assumptions \ref{assumption rational} and \ref{assumption order},
we solve the functional equation that corresponds to the stochastic recursion in terms of
survival function LSTs (Section \ref{Analysis of the stochastic recursion}).

\section{Duality}\label{section duality}

It turns out that it is more fruitful to regard the survival function defined in the previous section as a distribution rather than as a function. In this section we point out that the finite horizon survival functions $F_n^s$ can be seen as the c.d.f.s of the so called reflected version of the random walk. This reflected version satisfies a recursive equation that is then exploited in the subsequent sections in order to calculate the infinite horizon survival function.

For a one-dimensional random walk it is well known that its running maximum  has the same distribution as the reflected version of the random walk:

\begin{equation}\label{scalar inf sup}\sup\limits_{0\leq i\leq n} S_i^{(1)} \deq S_n^{(1)} - \inf\limits_{0\leq i\leq n} S_i^{(1)},\end{equation}
where $\deq$ denotes equality in distribution.
Here we could, e.g., take $S_n^{(1)} = - \sum_{i=1}^n X_i^{(1)}$, where $X_i^{(1)} =  A_i - B_i^{(1)}/c^{(1)}$
as defined in (\ref{RtnXn}) for the one-dimensional first risk reserve process.
In this section we show that this relation is still valid for the embedded random walk related to the general
{\em bivariate} process $\{-R_t\}_{t \geq 0}$, without having enforced the assumptions from the previous section. This relation is further used to derive a recursion for the two-dimensional transform of the survival function. We also give an interpretation of this recursion in terms of excursions away from the running maximum of the reserve process $\{R_t\}_{t \geq 0}$.

\vspace{0.2cm}
Throughout this section we will work with the natural order on $\mathbb R^{2}$. For $x:=(x^{(1)},x^{(2)}),y:=(y^{(1)},y^{(2)})\in \mathbb R^{2}$, denote

\[x\preceq y \Leftrightarrow x^{(1)}\leq y^{(1)},\,x^{(2)}\leq y^{(2)}.\] To keep notations short, set \[x\vee y := (\max(x^{(1)},y^{(1)}),\, \max(x^{(2)},y^{(2)})),\]

\[x\wedge y := (\min(x^{(1)},y^{(1)}),\, \min(x^{(2)},y^{(2)})).\]
Denote by $S_n:= \sum_{i=1}^n (-X_i)$, $n\geq 1$ and $S_0=0$, the origin of $\mathbb R^2$.
$X_n$ was  defined in (\ref{RtnXn}).
Ruin can only occur at arrival epochs, and since arrivals are simultaneous, we can introduce the maximum aggregate loss up to the $n^{th}$ arrival epoch $M_n:= \bigvee_{i=0}^n S_i$, so that we have the following relation for $\tau_\wedge$, the scaled version of the exit time that was defined in (\ref{naive time to ruin}):

\begin{equation}\label{maximum aggregate loss}  \{\tau_\wedge(u^{(1)},u^{(2)})>t_n \} = \{M_n \preceq (u^{(1)},u^{(2)})\},  \end{equation}
for $t_n$ the $n^{th}$ arrival epoch. Note also that $\tau_\wedge$
 can now be rewritten in terms of the order relation '$\succeq$':

\[\tau_\wedge(u^{(1)},u^{(2)}) = \inf \{t_n\,;R_{t_n} \nsucceq 0\,|\, R_0=(u^{(1)},u^{(2)})\}.\]

\vspace{0.2cm}
\noindent
 The duality approach described in the following, allows one to obtain a recursion for the two-dimensional survival function by turning this 'exit problem' for the embedded random walk into a so-called 'reflection problem' (Theorem \ref{remark duality}). This recursion is then solved in Section \ref{Analysis of the stochastic recursion} under Assumptions \ref{assumption rational} and \ref{assumption order}.
The behaviour of the bivariate reserve process $\{R_t\}_{t \geq 0}$ is similar to the behavior of the reserve process studied in \cite{APP1}, one of the main differences being that in our set-up this is not a Markov process anymore.


\vspace{0.2cm}
For fixed $n$  let $S^*_k:=S_n-S_{n-k}$, $k=0,...,n$, so that $S^*_0=0$, $S^*_1 = -X_n$, $S^*_2=-X_n-X_{n-1}$, etc. Thus $S_i^*$, $i\leq n$, is obtained from $S_i$ by circularly permuting its increments.
The following lemma is the 2-dimensional version of (\ref{scalar inf sup}). It remains valid in any number of dimensions by making some straightforward modifications.
\begin{Lemma}\label{Lemma duality}
For all $n\geq 0$,
\[ M_n \deq S_n - \bigwedge\limits_{i=0}^n S_i.\]
\end{Lemma}

\proof We can write

\[S_n- \bigwedge\limits_{i=0}^n S_i =S_n + \bigvee\limits_{i=0}^n (-S_i)= \bigvee\limits_{i=0}^n (S_n-S_i)= \bigvee\limits_{i=0}^n S^*_{n-i}. \]
Here we used the relation: $-(x\wedge y) = (-x)\vee (-y)$, $x,y\in \mathbb R^2$. The above are all sample-path identities. The final step is to remark that the joint distribution of $(S^*_1,...,S_n^*)$ is the same as the distribution of $(S_1,...,S_n)$ because the increments $X_i$ are i.i.d., and the proof is complete.

\vspace{0.2cm}

Denote $W_n:= S_n - \bigwedge\limits_{i=0}^n S_i$, $n\geq 0$, the random walk reflected inside the nonnegative quadrant. An important step in our analysis of the survival function is the following recursion for the  sequence $(W_n)_{n\geq 0}.$

\begin{Proposition}\label{Prop Lindley recursion} The sequence $(W_n)_{n\geq 0}$ satisfies the following recursion path-wise:
 \[ W_{n+1} = \left(W_n - X_{n+1}\right) \vee 0, \]
and initial condition $W_0=0$.
\end{Proposition}

Roughly speaking, as soon as one of the components of $S_n$ reaches a new minimum, the running infimum is updated accordingly and therefore the corresponding component of $W_n$ is set to zero.
\proof The proof follows by exploring all four possibilities, depending on the position of $W_n - X_{n+1}= S_{n+1} - \bigwedge \limits_{i=0}^{n} S_i $ relative to the origin.
For example, if $W_n - X_{n+1}$ is in the second quadrant, that is, if $S_{n+1}^{(1)} \leq \inf\limits_{i\leq n} S_i^{(1)}$ and $S_{n+1}^{(2)} \geq \inf\limits_{i\leq n} S_i^{(2)}$, then

\[    S_{n+1}^{(1)} = \inf\limits_{i\leq n+1} S^{(1)}_i, \mbox{ and }  S^{(2)}_{n+1}\geq \inf\limits_{i\leq n+1}S_i^{(2)}. \]
Therefore $W_{n+1} = S_{n+1} - \bigwedge\limits_{i=0}^{n+1} S_i = (0,\, S_{n+1}^{(2)} - \inf\limits_{i\leq n+1}S_i^{(2)})$, and remark that this is the same as $\left(W_{n} - X_{n+1}\right)\vee 0.$
The other cases follow by analogous considerations, which completes the proof.

\vspace{0.2cm}
 We can regard the finite horizon survival function

 \[F^s_n(u^{(1)},u^{(2)}):= \prob(\tau_\wedge(u^{(1)},u^{(2)})>t_n)\]
as the c.d.f. of a {\em survival measure}. 
Relation (\ref{maximum aggregate loss}), Lemma \ref{Lemma duality} and Proposition \ref{Prop Lindley recursion} imply that this survival measure is nothing else but the distribution  of the {\em reflected random walk} $W_n$ inside the non-negative quadrant of $\mathbb R^2$.
 \begin{Theorem}[Duality]\label{remark duality} The following identity relates the distribution of the reflected version of the random walk to the finite horizon survival functions of the reserve process:

\begin{equation}\label{2-dimensional duality} \prob(R_{t_i} \succeq 0,\,i=1,...,n\,|\,R_0=(u^{(1)},u^{(2)})) = \prob (W_n\preceq (u^{(1)},u^{(2)})\,|\,W_0=0)  .  \end{equation}

\end{Theorem}

\proof That $W_n$ is the reflected version of the random walk follows directly from the fact that it is the solution of the recursive equation in Proposition \ref{Prop Lindley recursion}. In view of Lemma \ref{Lemma duality} and (\ref{maximum aggregate loss}), the other statement is obvious, so this concludes the proof.

\begin{Remark}\label{queueing remark}
In fact $\{W_n\}_{n\geq 0}$ is the waiting time process in an initially empty two-dimensional queueing model with
two servers, each with its own queue, and with simultaneous arrivals and the same
$\{A_n,B_n^{(1)},B_n^{(2)}\}_{n \geq 0}$ input process as for the risk reserve process;
$A_n$ indicates interarrival time and $(B_n^{(1)},B_n^{(2)})$ denotes the vector of service requirements for the two servers.

One-dimensional instances of the above 'duality relation' are well known in the risk insurance/queueing literature, see for example \cite{AsmussenRuin}, p. 45, p. 161 etc.
\end{Remark}

\section{A functional equation}
\label{functional}
In this section we consider
the Laplace-Stieltjes transform of the survival function $ F_n^s(u^{(1)},u^{(2)})$; Theorem \ref{remark duality} and (\ref{2-dimensional duality}) imply that this equals the transform of the bivariate waiting time for the $n^{th}$ customer:
\begin{equation} \label{Laplace Stieltjes connection}\mean e^{-s_1 W^{(1)}_n - s_2 W_n^{(2)}}= \int e^{-s_1 u^{(1)}-s_2 u^{(2)}} {\rm d} F_n^s(u^{(1)},u^{(2)}),\,\,\,\, \mathcal Re\,s_i\geq 0,\,\,\,i=1,2.\end{equation}
Our main goal in this section is to obtain a recursion between the LSTs of $(W_{n+1}^{(1)},W_{n+1}^{(2)})$ and $(W_n^{(1)},W_n^{(2)})$ using the path-wise recursion in Proposition \ref{Prop Lindley recursion}.
But first we point out a sample-path relation between the reserve process and the reflected random walk which will be useful in the next section (the relation in Lemma \ref{Lemma duality} holds in distribution only). The event that any of the reserves is running at the maximum is the same as the event that the corresponding component of the reflected random walk is at 0.

\begin{Lemma}\label{lemma max-reflection}
\[ \{ W^{(i)}_{n+1} = 0\} = \{ R^{(i)}_{t_{n+1}} = u^{(i)} + \max(0,-S^{(i)}_1,...,-S^{(i)}_{n+1})\},\,\,\,\,\, i=1,2.  \]
\end{Lemma}

\proof $R^{(i)}_{t_{n+1}}= u^{(i)} - S^{(i)}_{n+1}$ and notice the following equivalence:

\[  S^{(i)}_{n+1} - \min(0,S^{(i)}_1,...,S_n^{(i)},S^{(i)}_{n+1}) = 0 \Leftrightarrow -S^{(i)}_{n+1} = \max(0,-S^{(i)}_1,...,-S_n^{(i)},-S^{(i)}_{n+1}). \]

\begin{Remark} The event $\{W^{(i)}_n=0\}$ does not depend on the initial capital. The event on the RHS in Lemma \ref{lemma max-reflection} above does not restrict the reserve process to staying above 0; equivalently $W^{(i)}_n$ on the LHS is not restricted to staying below level $u^{(i)}$. This is in line with (\ref{2-dimensional duality}) in Theorem \ref{remark duality}.
\end{Remark}

Below we point out how one can obtain a recursive equation for the LST of the survival function $F_n^s(u^{(1)},u^{(2)})$ in the general case without the ordering assumption. In Section \ref{Analysis of the stochastic recursion} we solve this equation for the case when the risks are ordered.
\newline The bivariate recursion in Proposition \ref{Prop Lindley recursion} becomes in terms of LSTs:

\begin{align*}  \mean e^{-s_1 W_{n+1}^{(1)}-s_2 W_{n+1}^{(2)}} = \mean e^{-s_1(W^{(1)}_n -X^{(1)}_{n+1})^+ - s_2(W^{(2)}_n -X^{(2)}_{n+1})^+ },
 \end{align*}
 with $(x)^+$ denoting $\max(x,0)$. Hence, with $1_{\{E\}}$ the indicator function of event $E$,

  \begin{align}\label{general functional equation} \mean e^{-s_1 W_{n+1}^{(1)}-s_2 W_{n+1}^{(2)}} = &\mean \left[e^{-s_1(W^{(1)}_n -X^{(1)}_{n+1}) - s_2(W^{(2)}_n -X^{(2)}_{n+1})}1_{\{X^{(1)}_{n+1}< W^{(1)}_n,\,X^{(2)}_{n+1}< W^{(2)}_n\}}\right] \notag \\
   &+ \mean \left[e^{-s_1(W^{(1)}_n -X^{(1)}_{n+1})}1_{\{X^{(1)}_{n+1}< W^{(1)}_n,\, X^{(2)}_{n+1} \geq W^{(2)}_n\}}\right] \notag\\
   &+ \mean \left[e^{-s_2(W^{(2)}_n -X^{(2)}_{n+1})}1_{\{X^{(1)}_{n+1}\geq W^{(1)}_n,\, X^{(2)}_{n+1} < W^{(2)}_n\}} \right] \notag\\
   &+ \prob(X^{(1)}_{n+1} \geq W^{(1)}_n , X_{n+1}^{(2)} \geq W_{n}^{(2)}).
  \end{align}



In view of (\ref{Laplace Stieltjes connection}), the left-hand side of
(\ref{general functional equation}) represents the LST
of the survival measure $F^s_{n+1}$.
Below we also interpret each of the four terms in the righthand side in terms of transforms of survival measures.
\\
{\em Term 1:}
In terms of excursions away from the maximum reserve (Lemma \ref{lemma max-reflection}),
the first term on the RHS represents the transform of the survival measure in the event that both reserves are during an excursion below the running maximum at time $t_{n+1}$,

\[   \mean e^{-s_1 W^{(1)}_{n+1} - s_2 W^{(2)}_{n+1} }1_{\{W^{(1)}_{n+1}>0,\, W^{(2)}_{n+1}>0\}} =  \int e^{-s_1 u^{(1)}-s_2 u^{(2)}} 1_{\{R_{t_{n+1}} \prec \bigvee\limits_{k=0,...,n} R_{t_k}\}} {\rm d} F_{n+1}^s(u^{(1)},u^{(2)}). \]
Above we used that on the event $\{W^{(i)}_{n+1}>0\}$, it holds that $W^{(i)}_{n}- X^{(i)}_{n+1}=W^{(i)}_{n+1}$ (Proposition \ref{Prop Lindley recursion}).
\\
{\em Terms 2 and 3}: These can be translated in terms of survival functions using (\ref{Laplace Stieltjes connection}) again:

\[\mean e^{-s_1 W^{(1)}_{n+1}}1_{\{ W^{(1)}_{n+1}>0,\, W^{(2)}_{n+1}=0\}} = \int\limits_{0+}^\infty e^{-s_1 u^{(1)}}1_{\{R^{(1)}_{t_{n+1}} < \max\limits_{k=0,...,n} R^{(1)}_{t_k};\,R^{(2)}_{t_{n+1}} \geq \max\limits_{k=0,...,n} R^{(2)}_{t_k} \}}  {\rm d}F^s_{n+1}(u^{(1)},0),\]

\[\mean e^{-s_2 W^{(2)}_{n+1}}1_{\{W^{(1)}_{n+1}=0,\,W^{(2)}_{n+1}> 0\}} = \int\limits_{0+}^\infty e^{-s_2 u^{(2)}} 1_{\{R^{(1)}_{t_{n+1}} \geq \max\limits_{k=0,...,n} R^{(1)}_{t_k};\,R^{(2)}_{t_{n+1}} < \max\limits_{k=0,...,n}R^{(2)}_{t_k} \}} {\rm d}F^s_{n+1}(0,u^{(2)}), \]
are 'boundary' transforms. 
The survival measure that corresponds to $F^s_{n}$ can have positive mass on the axes of the non-negative quadrant, if nowhere else, at least $F_{n}^s(0,0)=\prob(W_n=0\,|\,W_0=0)$ is positive, i.e.,
it has an atom in the origin.
\\
{\em Term 4}:
This is the probability that both reserves are surviving and running at a maximum, which by Lemma 2 and Theorem 1 is $F_{n+1}^s(0,0)$.


\vspace{0.2cm}
 From an analytic point of view it is more convenient to rewrite (\ref{general functional equation}) as a recursion. After adding and subtracting appropriate terms, one obtains, 

\begin{align}\label{general functional equation 2} \mean e^{-s_1 W_{n+1}^{(1)}-s_2 W_{n+1}^{(2)}} &=\mean e^{s_1X^{(1)}_{n+1} + s_2X^{(2)}_{n+1}}\, \mean e^{-s_1W^{(1)}_n  - s_2W^{(2)}_n } \notag \\
   &+ \mean \left\{e^{-s_1(W^{(1)}_n -X^{(1)}_{n+1})}[1- e^{-s_2(W^{(2)}_n - X^{(2)}_{n+1})}]1_{\{W^{(1)}_{n+1}>0,\,W^{(2)}_{n+1}=0\}} \right\} \notag\\
   &+ \mean \left\{ e^{-s_2(W^{(2)}_n -X^{(2)}_{n+1})}[1- e^{-s_1(W^{(1)}_n - X^{(1)}_{n+1})}]1_{\{W^{(1)}_{n+1}=0,\,W^{(2)}_{n+1}>0\}}\right\}  \notag\\
   &+ \mean \left\{ [1- e^{-s_1(W^{(1)}_n - X^{(1)}_{n+1})-s_2(W^{(2)}_n - X^{(2)}_{n+1})}]1_{\{W^{(1)}_{n+1}=0,\,W_{n+1}^{(2)}=0\}}\right\}.
  \end{align}
 Above we used that the increment $X_{n+1}$ is independent of $W_n$.
Under the assumption that the vector $W_n$ has a limit in distribution, $W$, as $n\rightarrow \infty$, (\ref{general functional equation 2}) becomes

\begin{align}\label{queueing friendly formula}
K(s_1,s_2) \mean e^{-s_1 W^{(1)}-s_2 W^{(2)}} &=\mean \left\{ e^{-s_1(W^{(1)}-X^{(1)})}[1- e^{-s_2(W^{(2)} - X^{(2)})}]1_{\{W^{(1)}>X^{(1)},\,W^{(2)}\leq X^{(2)}\}}\right\} \notag\\
   &+ \mean \left\{ e^{-s_2(W^{(2)} -X^{(2)})}[1- e^{-s_1(W^{(1)} - X^{(1)})}]1_{\{W^{(1)}\leq X^{(1)},\,W^{(2)}>X^{(2)}\}}\right\}  \notag\\
   &+ \mean \left\{ [1- e^{-s_1(W^{(1)} - X^{(1)})-s_2(W^{(2)} - X^{(2)})}]1_{\{W^{(1)} \leq X^{(1)},\,W^{(2)}\leq X^{(2)}\}} \right\},
\end{align}
 with "kernel" $K(s_1,s_2):= 1- \mean e^{s_1X^{(1)}+s_2X^{(2)}}$ and $\mathcal Re\,s_i=0$, $i=1,2$.  

\begin{Remark}\label{weak limit}
From Lemma \ref{Lemma duality}, $W$ has the same distribution as the all-time supremum $M:=\lim\limits_{n\rightarrow \infty} M_n$; and $M_n$ being a sequence of nondecreasing random vectors  w.r.t. the order '$\preceq$', the limit always exists a.s., although it may have a defective distribution.
In the next section we give a sufficient condition for $M$ to have a proper distribution under the assumption that risks are ordered.
\end{Remark}

\section{Wiener-Hopf analysis of the stochastic recursion}\label{Analysis of the stochastic recursion}

In this section we resolve the functional equation (\ref{queueing friendly formula}) under the Assumptions \ref{assumption rational} and \ref{assumption order}, i.e., we find the LST of the infinite horizon survival function:
\[ F^s= \lim\limits_{n\rightarrow \infty} F^s_n,\]
the limit being considered in distribution. Theorem  \ref{remark duality} together with a limit argument shows that this weak limit is the same as the c.d.f. of the stationary version of the waiting time process $(W_n)_{n\geq 0}$ (see also Remark \ref{queueing remark}).

The section is divided into three subsections.
Subsection \ref{subs5.1} prepares the ground, by making a key observation about the functional equation
(\ref{queueing friendly formula}), introducing some notation and discussing the stability condition.
Subsection \ref{subs5.2} expresses the two-dimensional LST $\psi(s_1,s_2)$ of $F^s$ in a one-dimensional unknown function $C(s_1)$ (Proposition
\ref{Proposition W-H 1 parameter}).
That function is determined in Subsection \ref{subs5.3}, yielding our main result: Theorem \ref{Main result}.

\subsection{Preparations}
\label{subs5.1}
Introduce the extra claim amount $D_n:=B_n^{(1)}/c^{(1)} - B_n^{(2)}/c^{(2)}= X_n^{(2)} - X_n^{(1)}$,  so that  the increments of the random walk $S_n$ can be represented as $-X_n =(-X_n^{(2)}+D_n,-X_n^{(2)})$.

We first make the following {\em key observation}:
The ordering assumption (\ref{ordering assumption}) implies that when $R^{(1)}_{t_n}$ is at a maximum, necessarily $R^{(2)}_{t_n}$ is at a maximum. Via Lemma \ref{lemma max-reflection}, this corresponds to the fact that the events
 \[\{W^{(1)}_{n}\leq  X_{n+1}^{(1)},\,W^{(2)}_{n} > X_{n+1}^{(2)}\} = \{W^{(1)}_{n+1}=0,W^{(2)}_{n+1}>0\}\]
 are null for all $n\geq 0$. This means that the third term on the RHS of (\ref{general functional equation 2}) is null and hence the second term on the RHS of (\ref{queueing friendly formula}) vanishes as well, so that after regrouping terms, (\ref{queueing friendly formula}) can be rewritten as

\begin{equation}\label{functional identity1}  K(s_1,s_2) \psi(s_1,s_2) = - \psi_{1} (s_1,s_2) + \psi_2(s_1) + {\textstyle\prob(X^{(1)}\geq W^{(1)})},  \end{equation}
 where
\begin{eqnarray*}  \psi_{1} (s_1,s_2) &= &\mean e^{-s_1(W^{(1)}-X^{(2)}  + D) - s_2(W^{(2)} -X^{(2)})}1_{\{X^{(2)}\geq W^{(2)} \}},  \\
  \psi_2(s_1) &= &\mean e^{-s_1(W^{(1)}-X^{(2)}  + D)} 1_{\{W^{(1)}+D > X^{(2)}\geq W^{(2)}\}}.   \end{eqnarray*}

\vspace{0.2cm}

Consider the following function:
\[ \tilde K(s_1,z):= 1 - \mean e^{-s_1 D + z X^{(2)}},\;\;\; \mathcal Re\;s_1\geq 0,\;\mathcal Re\;z=0. \]
This is related to $K(s_1,s_2)$ that appears in (\ref{queueing friendly formula}) through a change of coordinates: $\tilde K(s_1,z)=K(s_1,z-s_1)$. In addition, remark that for fixed $z$, $\tilde K(s_1,z)$ is indeed analytic in $\mathcal Re\, s_1 >0$ because $D\geq 0$ a.s.
Now let us change the coordinates: $(s_1,s_2)\rightarrow (s_1,s_1+s_2)=:(s_1,z),$ and denote $\tilde \psi(s_1,z):=\psi(s_1,z-s_1)$. Then $\psi_{1}(s_1,s_2)$ becomes

\[ \psi_{1}(s_1,s_2) = \mean \left[e^{-s_1(W^{(1)}- W^{(2)} +D ) - z(W^{(2)}-X^{(2)} )}1_{\{X^{(2)}\geq W^{(2)} \}}\right] =: \tilde\psi_{1}(s_1,z), \]

\noindent and therefore (\ref{functional identity1}) can be rewritten as

\begin{equation}\label{functional identity1.5} \tilde K(s_1,z) \tilde\psi(s_1,z) = - \tilde\psi_{1} (s_1,z) + \psi_2(s_1) + \prob(X^{(1)}\geq W^{(1)}).  \end{equation}
\paragraph{\textbf{Running example:}} One of the simplest examples that we will use throughout is obtained when taking the joint distribution  of $(A,D,B^{(2)})$ to be such that conditional on a random variable $N$, these $A$, $D$ and $B^{(2)}$ are independent and have Erlang distributions of order $N$ and rates respectively $\lambda$, $\mu_D$ and $\mu$.
To keep things as simple as possible, we choose $\prob(N=1) = \prob(N=2) = \frac{1}{2}$ and rates $\lambda=1 $, $\mu_D=3$, $\mu=2$; we also choose the income rates $c^{(1)}$ and $c^{(2)}$ to be equal to 1.

The kernel $\tilde K(s_1,z)$  has the following simple form:

\begin{equation}\label{running example} \tilde K(s_1,z)  = 1- \frac{3(3+s_1)(1-z)(2+z) + 18}{(3+s_1)^2(1-z)^2(2+z)^2}. \end{equation}
From a specific example like this it becomes clear that the coefficients of the numerator, for example, as a polynomial in $z$ are themselves polynomials in $s_1$ and vice versa.

\vspace{0.2cm}
We are now ready to formulate a Wiener-Hopf boundary value problem in variable $z$.
For fixed $s_1$, $\mathcal Re\; s_1>0$, $\tilde \psi_{1}(s_1,z)$ is analytic in $\mathcal Re\; z <0$ (by analytic continuation), while $\tilde \psi(s_1,z)$ is analytic (by analytic continuation) in $\mathcal Re\; z>0$. These statements follow easily from the probabilistic nature of these functions.
In particular, notice that
$\tilde \psi(s_1,z) = \mean {\rm e}^{-s_1(W^{(1)} - W^{(2)}) -zW^{(2)}}$.
Also on the event $\{W^{(2)}<X^{(2)}\}$, the random variable $e^{-z(W^{(2)}-X^{(2)})}$ is uniformly bounded in $\mathcal Re\,z\leq 0$, hence the analyticity of $\tilde \psi_1(s_1,z)$ follows by an application of Lebesgue's dominated convergence theorem.

\vspace{0.2cm}
The approach we take in order to solve (\ref{functional identity1.5}) uses a Wiener-Hopf factorization with a parameter. More precisely, for each fixed $s_1$, $\mathcal Re\,s_1>0$, we will factorize the bivariate kernel $\tilde K(s_1,z)=\tilde K_{s_1}(z)$ that appears in (\ref{functional identity1.5}) into $\tilde K_{s_1}(z) = \tilde K^+_{s_1}(z) \tilde K^-_{s_1}(z)$, such that $\tilde K_{s_1}^+(z)$ can be  analytically continued in $\mathcal Re\, z>0$ and $\tilde K^-_{s_1}(z)$ can be analytically continued in $\mathcal Re\,z <0$.
The Wiener-Hopf factorization that solves (\ref{functional identity1.5}) is discussed in the next subsection.
Finally, once we resolve (\ref{functional identity1.5}), the solution to (\ref{functional identity1}) follows by reverting to the original coordinate system $(s_1,s_2)$.

\begin{Remark}\label{remark appology}
A reason why we prefer the notation using the argument $s_1$ as a subscript $\tilde K^{\pm}_{s_1}(z)$ for these factors is that they are in general obtained by pasting together different branches of multi-valued complex functions in $s_1$ using analytic continuation. More precisely, since $\tilde K(s_1,z)$ is a rational function, the 1-parameter Wiener-Hopf factors may have branch cuts in $\mathcal Re\,s_1>0$ (discontinuities) as functions of the argument $s_1$; then as it follows from Proposition \ref{Proposition W-H 1 parameter} below, we must choose the values of the zeroes of the kernel that have positive real part for $\mathcal Re\,s_1>0$ and glue them together (using analytic continuation). Because of this, the 1-parameter Wiener-Hopf factors $\tilde K^+_{s_1}(z)$ and $\tilde K^-_{s_1}(z)$ are not functions of $s_1$ in the real sense. This will be the case with the zeroes of the kernel from Example 2 in Section \ref{Section numerics}.
\end{Remark}

 \vspace{0.2cm}

Finally a word about conditions under which the limiting distribution of the two-dimensional waiting time process $\{W_n^{(1)},W_n^{(2)}\}_{n=1,2,\dots}$ exists,
or equivalently, under which survival of both risk reserves has a positive probability.
It will turn out from the analysis below that a necessary condition for the existence of a proper limit in distribution $W$ is $\rho_1:=\mean B^{(1)}/(c^{(1)}\mean A)<1$. This is easy to interpret in our case, because it is sufficient to ensure positive safety loading for line 1 which receives always larger claims - we then automatically have positive safety loading for the second insurance line.
 The safety loading condition for the second reserve process $\rho_2:=\mean B^{(2)}/(c^{(2)}\mean A)<1$ will be necessary for the Wiener-Hopf factorization to hold.
 The two Wiener-Hopf factors will be initially determined up to a certain unknown 'boundary' function $C(s_1)$ that appears in Equation (\ref{intermed Liouville}) below; we further determine this boundary function by noting that the marginal reserve process $R_t^{(1)}$ behaves as a (one-dimensional) generalized Sparre-Andersen reserve process with dependence between inter-arrival times and subsequent claim sizes, for which an analysis of the survival function is available in \cite{paper1}. At this point the safety loading condition $\rho_1<1$ becomes necessary.

\subsection{A Wiener-Hopf factorization}
\label{subs5.2}
In this subsection we determine the double transform $\tilde \psi(s_1,z)$ up to a -- yet -- unknown one-dimensional function $C(s_1)$.
In the next subsection we will determine this function, which turns out to be related to the first insurance line only.
\begin{Proposition}[\textbf{Wiener-Hopf factorization with a parameter}]\label{Proposition W-H 1 parameter} Under  Assumption \ref{assumption rational} the double LST $\tilde \psi(s_1,z)$ is of the form

\begin{equation} \label{intermed Liouville} \tilde \psi(s_1,z) = \mean e^{-s_1(W^{(1)}- W^{(2)}) -z W^{(2)}} = C(s_1)\, \tilde K^+_{s_1}(z)^{-1}. \end{equation}
$C(s_1)$ is a yet to be determined analytic function, $\mathcal Re\, s_1>0$. For fixed $\mathcal Re\, s_1>0$, $\tilde K^+_{s_1}(z)$ is analytic for $\mathcal Re\, z>0$, continuous up to the boundary and it factorizes  $\tilde K(s_1,z)$ into

\[ \tilde K(s_1,z) = \tilde K^+_{s_1}(z) \,\tilde K^-_{s_1}(z),\]
such that $\tilde K^-_{s_1}(z)$ is analytic for $\mathcal Re\, z<0$ and continuous up to the boundary.
\end{Proposition}

\proof Under the above Assumption \ref{assumption rational}, the triple transform of $A$, $B^{(1)}/c^{(1)}$, $B^{(2)}/c^{(2)}$ is \newline $H(q_0,q_1/c^{(1)}, q_2/c^{(2)})$ and the kernel is $\tilde K(s_1,z) = 1- H(-z,s_1/c^{(1)},(z-s_1)/c^{(2)})$  which is also a rational function with  representation  $\tilde K(s_1,z):=1- \frac{f(s_1,z)}{g(s_1,z)}$. The assumption that $A>0$ a.s. implies that the degree of $f_{s_1}(z)$ is strictly less than the degree of $g_{s_1}(z)$ as polynomial functions in the argument $z$ (see the discussion in Remark \ref{remark appology}).  Now the functional equation  (\ref{functional identity1.5}) becomes

\begin{equation}\label{functional identity2}  \frac{g(s_1,z)- f(s_1,z)}{g(s_1,z)} \tilde\psi(s_1,z) = - \tilde\psi_{1} (s_1,z) + \psi_2(s_1) + \prob(X^{(1)} \geq W^{(1)}).  \end{equation}

The first step is to factorize the kernel into two factors with respect to the $z$ variable and regroup (\ref{functional identity2}) into an analytic function in $\mathcal Re\; z>0$ on the LHS and an analytic function in $\mathcal Re\; z<0$ on the RHS. Once we have this representation, we can use Liouville's Theorem to determine both sides of (\ref{functional identity2}) up to the function $C(s_1)$ .

 Remove all the poles with nonnegative real part from the LHS of (\ref{functional identity2}).
 We keep for now $s_1$ fixed with $\mathcal Re\;s_1\geq0$, and denote by

 \[ g_{s_1}^-(z):=\prod_{i: Re\;z_i(s_1)<0}(z -z_i(s_1)),\]
 where $z_i(s_1)$ are zeroes of $g(s_1,z) = g_{s_1}(z)$; also put $g_{s_1}^+(z):= \frac{g_{s_1}(z)}{g_{s_1}^-(z)}$, so that we have the factorization $g_{s_1}(z) = g_{s_1}^+(z)\, g_{s_1}^-(z)$.
Upon multiplying both sides of (\ref{functional identity2}) by $g_{s_1}^+(z)$, the LHS becomes analytic for all $\mathcal Re\,z>0$ and continuous up to the imaginary axis.
 Similarly, the RHS is analytic for all $\mathcal Re\,z<0$ and continuous for $\mathcal Re\,z\leq 0$.
 Since these two coincide for $\mathcal Re\;z=0$, they are analytic continuations of each other, in particular $\frac{g_{s_1}(z)- f_{s_1}(z)}{g^-_{s_1}(z)} \tilde \psi_{s_1}(z)$ is an entire function in $z$.
 Because $\deg f_{s_1}(z)\leq \deg g_{s_1}(z)$, asymptotically
\[ \frac{g_{s_1}(z)- f_{s_1}(z)}{g^-_{s_1}(z)} \tilde \psi_{s_1}(z)  = O(z^{m_+(s_1)}),
\]
where $m_+(s_1):=\deg g_{s_1}^+(z)$.
By virtue of Liouville's theorem (\cite{Titchmarsh}, p.85),

\begin{equation}\label{Liouville identity}\tilde \psi_{s_1}(z)  =  \frac{g^-_{s_1}(z)}{g_{s_1}(z)- f_{s_1}(z)}  P_{s_1}(z), \end{equation}
where (for fixed $s_1\geq 0$), $P_{s_1}(z)$ is a polynomial in $z$ with $\deg P_{s_1}(z)\leq m_+(s_1).$ From (\ref{Liouville identity}) it follows immediately that $P_{s_1}(z)$ must have all the zeroes with non-negative real part of the denominator $g_{s_1}(z)- f_{s_1}(z)$. Now a key part in the argument is the fact that the denominator $g_{s_1}(z)-f_{s_1}(z)$ has the same number of zeroes in $\mathcal Re\,z\geq 0$ as $g^+_{s_1}(z)$. The proof of this statement is deferred to the Appendix in Proposition \ref{Rouche problem}. Thus we have $\deg P_{s_1}(z)\geq m_+$. Together with the upper bound on the degree of $P_{s_1}(z)$, this implies $\deg P_{s_1}(z)=m_+$; moreover, this determines $P_{s_1}(z)$ up to a constant factor (constant being relative to $z$ !)

\[ P_{s_1}(z) = C(s_1) \prod_{i: \mathcal Re\; v_i(s_1)\geq 0} (z-v_i(s_1)), \] where $v_i(s_1)$ are zeroes of $g_{s_1}(z)- f_{s_1}(z)$. Upon replacing the above in (\ref{Liouville identity}), we have found the one-parameter positive Wiener-Hopf factor

\begin{equation}\label{positive W-H factor} \tilde K_{s_1}^+(z) = \frac{\prod\limits_{i: \mathcal Re\; v_i(s_1) < 0} (z-v_i(s_1))}{g_{s_1}^-(z)}.  \end{equation}
And in particular the above also determines $\tilde K_{s_1}^-(z)$:

\[ \tilde K_{s_1}^-(z)= \frac{\tilde K(s_1,z)}{\tilde K^+_{s_1}(z)}, \]
and the proof is complete.

\vspace{0.2cm}
\paragraph{\textbf{Running example:}} For fixed $s_1$, we can carry out the factorization for the kernel (\ref{running example}) in the running example. Below we give the zeroes of the numerator

 \[ v_1(s_1)=\frac{-s_1-3-\sqrt{3}\sqrt{(s_1+3)(1+3s_1)}}{2(3+s_1)},\, v_2(s_1)=\frac{-s_1-3+\sqrt{3}\sqrt{(s_1+3)(3s_1+1)}}{2(3+s_1)},\]
 \[   v_3(s_1)= \frac{-s_1-3-\sqrt{(s_1+3)(3s_1+13)}}{2(3+s_1)},\, v_4(s_1)=\frac{-s_1-3+\sqrt{(s_1+3)(3s_1+13)}}{2(3+s_1)} .\]

The radicals above are defined when the cut in the complex plane is taken along the negative half of the real axis and the complex arguments are measured from $-\pi$ to $\pi$. This convention determines the so-called principal value of the square root function. The negative real half-axis will be a discontinuity line for the square root function $\sqrt{z}$ as a function of a complex variable because as the argument $z$ approaches the negative real half-axis,

\[  (z^--z_0)^{1/2}= e^{i\pi}(z^+-z_0)^{1/2}, \]
where $z_0$ lies on the negative real half-axis, that is, $\mathcal Re\,z_0<0$ and $\mathcal Im\,z_0=0$;  $z^{\pm}$ denotes the limit of $z$ towards $z_0$ from respectively above and below the real axis. We will call such lines of discontinuity branch cuts.
The branch points of $v_1$ (and $v_2$) are $-3$ and $-1/3$. The branch cuts are then the curves generated by the equations
\[  \arg(s_1+3) + \arg(1+3s_1) = \pm \pi.\]

It is a problem of plane geometry to see that these branch cuts are constituted by the line segment joining the two branch points, together with the perpendicular line on this segment that passes through its mid-point. The situation is similar for $v_3(s_1)$ and $v_4(s_1)$.

  By inspecting the zeroes of the numerator for $\mathcal Re\,s_1>0$, exactly $v_1(s_1)$ and  $v_3(s_1)$ are negative, where by positive/negative values of complex numbers we will always refer to their real parts. There are as many negative zeroes in the denominator, which is already confirmed by Proposition \ref{Rouche problem}. 
  Moreover, the branch cuts of neither $v_1$ nor $v_3$ are located in the right half-plane, which means these zeroes are regular functions for positive values of $s_1$.

  Having isolated the negative zeroes, the one-parameter positive Wiener-Hopf factor from (\ref{positive W-H factor}) is

\begin{equation}\label{one-parameter WH example}\tilde K^+_{s_1}(z)= \frac{(z-\frac{-s_1-3-\sqrt{3}\sqrt{(s_1+3)(1+3s_1)}}{2(s_1+3)})(z-\frac{-s_1-3-\sqrt{(s_1+3)(3s_1+13)}}{2(s_1+3)} )}{(z+2)^2}.\end{equation}

\begin{Remark} Interestingly, $\tilde K^+_{s_1}$ is  not a rational function anymore in the argument $s_1$ (however, in this example it is meromorphic in the argument $s_1$, for $\mathcal Re\;s_1>0$).
A queueing theoretic explanation of this remark can be found by comparing the double transform (\ref{original transform}) obtained below with the decomposition results for a particular case of the present model in \cite{paper2}. For the  process with Markov arrivals studied therein, the stationary waiting time stochastically decomposes into two components, one of which is related to the extra busy period length in the longest queue and it is well known that
busy periods in general do not have rational transforms (already in the case of an $M$/$M$/1 system, the busy period has a non-rational transform).

\end{Remark}

Moreover, it is easy to see that the marginal factor $\tilde K_{s_1}^+(s_2)|_{s_1=0}$ is a rational function. From this and (\ref{original transform}) below follows that the marginal transform $\psi(0,s_2)$ is a rational function. Also, $\psi(s_1,0)$ is a rational function because for $s_2=0$ the factor $\tilde K^+_{s_1}$ cancels against itself in (\ref{original transform}).
The rationality of the marginal transforms  is clear from their queueing interpretation because these are the transforms of the univariate survival functions/waiting times for the two insurance lines/queueing systems in isolation (see the discussion in \cite{CohenSingleServer} p.325 and the references therein).

\subsection{The main result}
\label{subs5.3}
We are now ready to formulate and prove the main result.

\begin{Theorem}\label{Main result}
Under the safety loading condition for the riskier line 1, $\rho_1<1$, the infinite horizon survival function $F^s(u^{(1)},u^{(2)})$ is a (proper) probability distribution function with support the non-negative quadrant in $\mathbb R^{2}$, and its LST is given by


\begin{equation} \label{original transform} \psi(s_1,s_2) = \int e^{-s_1 u^{(1)}-s_2 u^{(2)}} {\rm d} F^s(u^{(1)},u^{(2)}) =  \frac{ K_{pr}^+(0)}{K^+_{pr}(s_1)}\; \frac{\tilde K^+_{s_1}(s_1)}{\tilde K^+_{s_1}(s_1+s_2)}, \end{equation}
for $\mathcal Re\,s_i\geq 0,\,i=1,2$. $\tilde K_{s_1}^+(z)$ is given in (\ref{positive W-H factor}) and here is evaluated at $z=s_1$ and at $z=s_1+s_2$. $K^+_{pr}(s_1)$ is the positive Wiener-Hopf factor of the projected one-dimensional kernel $\tilde K(s_1,s_1)= K(s_1,0)$, i.e., the unique function analytic in $\mathcal Re\, s_1>0$, continuous in $\mathcal Re\, s_1\geq 0$ that factorizes $K(s_1,0)$ into

\[ K(s_1,0) = K^{+}_{pr}(s_1)\,K^-_{pr}(s_1),  \]
with $K^-_{pr}(s_1)$ analytic in $\mathcal Re\, s_1<0$ and continuous in $\mathcal Re\, s_1\leq 0$ (see for instance \cite{Prabhu80}, Thm.7 p.55). Under Assumption \ref{assumption rational} it is of the form (see also~\cite{paper1})

\begin{equation}\label{proj factor} K^+_{pr}(s_1) = \frac{\prod_{j} (s_1 - \tilde v_j^-)}{\prod_{j}(s_1 - v_j^-)},\end{equation}
with $\tilde v_j^-$ the negative zeroes of $K(s_1,0) $ and $v_j^-$ its negative poles.

\end{Theorem}

\proof
Our starting-point is
(\ref{intermed Liouville}), and our goal is
to determine the one yet unknown function $C(s_1)$ in that formula. The idea is that, since $C(s_1)$ stays the same irrespective of the value of $z$, we are free to choose any $z$. Since by definition, $\tilde K^+_{s_1}(z)$ is analytic for all $\mathcal Re\, z> 0$ and continuous in $\mathcal Re\, z\geq 0$, take $z=s_1$:

\begin{equation}\label{marginal Liouville} \tilde \psi(s_1,s_1) = \mean e^{-s_1 W^{(1)}} = C(s_1) [\tilde K^+_{s_1}(s_1)]^{-1}.  \end{equation}

We can determine $C(s_1)$ from (\ref{marginal Liouville}), because $\tilde \psi(s_1,s_1) = \psi(s_1,0)=\mean e^{-s_1 W^{(1)}}$ is the steady-state waiting time transform in the marginal $G/G/1$ queue with dependence between inter-arrival times and service requirements, which has been determined in \cite{paper1}.
That paper was devoted to an analysis of a one-dimensional risk/queueing model that amounts to the present model with $B^{(2)} \equiv 0$.
The kernel of the functional identity for this marginal queue with generic service requirement $B^{(1)}$ and correlated inter-arrival time $A$ is $\tilde K(s_1,s_1)$; the corresponding Rouch\'e problem is to prove that $g(s_1,s_1)$ and $g(s_1,s_1)- f(s_1,s_1)$ have the same number of nonnegative zeroes.
This has been carried out in \cite{paper1}. Its Formula (6) reads

\begin{equation}\label{marginal waiting time}\tilde \psi(s_1,s_1) =\mean e^{-s_1W^{(1)}}=  K^+_{pr}(0)[ K^+_{pr}(s_1)]^{-1}, \end{equation}
with $K^+_{pr}(s_1)$ the positive Wiener-Hopf factor of the projected kernel $K(s_1,0)$:

\[ K^+_{pr}(s_1) = \frac{\prod_{j} (s_1 - \tilde v_j^-)}{\prod_{j}(s_1 - v_j^-)},\]
such that $\tilde v_j^-$ are the negative zeroes of $K(s_1,0) $ and $v_j^-$ its negative poles, and the normalizing constant $K^+_{pr}(0)$ is equal to the atom at 0 of $W^{(1)}$:

\[\prob(W^{(1)}=0)= \prod_{k} (- v^-_k) /\prod_{j} (- \tilde v^-_j).\]
Formula (\ref{marginal waiting time}) together with (\ref{marginal Liouville}) now determine $C(s_1)$:

\begin{equation}\label{C(s) formula}C(s_1)= K^+_{pr}(0) \, [ K^+_{pr}(s_1)]^{-1}\,\tilde K^+_{s_1}(s_1),\end{equation}
and with this we obtain the transform of the joint waiting time distribution, or equivalently
of the survival function $F^s(u^{(1)},u^{(2)})$, from (\ref{intermed Liouville}), upon switching back to the original coordinates. The proof is complete.

\begin{Remark} An important remark is that the factors in (\ref{original transform}), $\tilde K^+_{s_1}(s_1)$ and $K^+_{pr}(s_1)$, as defined in Proposition \ref{Proposition W-H 1 parameter} and in Theorem \ref{Main result}, are not the same. One can already compare (\ref{one-parameter WH example}) with (\ref{proj factor ex1}) for Example 1 in the following section.

\noindent More precisely, the operations of taking the projection and carrying out the factorization do not commute with each other, in contrast with the one-dimensional Fluctuation Theory of random walks. See also Section 13 in \cite{KingmanAlgebra}.

\noindent $\tilde K^+_{s_1}(s_1)$ is defined by first carrying out the one-parameter Wiener-Hopf factorization for $\tilde K(s_1,z)$ and then projecting the positive factor onto the main diagonal of the 2-dimensional complex space: $z=s_1$.

\noindent  On the other hand, $K^+_{pr}(s_1)$ is obtained by first projecting $\tilde K(s_1,z)$ onto the main diagonal $z=s_1$ and then carrying out the Wiener-Hopf factorization for the projected kernel $\tilde K(s_1,s_1)= K(s_1,0)$.
\end{Remark}

\begin{Remark} The LST in (\ref{original transform}) has a product form. It can be shown that the bivariate LST of the reflected random walk  decomposes in a similar way as in Theorem 1 of \cite{paper2},  where one of the factors is related to a modified workload process.
\end{Remark}

\section{Examples and numerical inversion}\label{Section numerics}

In the previous sections we dealt with some theoretical aspects related to obtaining the transform of the survival function. It turns out that additional insight can be obtained by applying the previous results to some specific examples. Our aims in this section are:

(i) to provide examples for which the Laplace-Stieltjes transforms of the survival measures can be calculated, based on the general results obtained in the previous section. 

(ii) to explain the various analytic challenges that appear when one tries to determine the Laplace-Stieltjes transform of the survival measure for some specific classes of distributions for the input $(A,D,B^{(2)})$.

(iii) numerical inversion of the bivariate Laplace-Stieltjes transform in (\ref{original transform}) and the comparison between the risks for various possible correlations between the claim sizes $(B^{(1)}_n,B^{(2)}_n)$ and inter-arrival times $A_n$.

We begin by explaining the inversion algorithm and how we applied it. However, in order to obtain the input for the algorithm, we need to follow the steps in Section 5 and construct the Wiener-Hopf factors. It turns out that this presents a challenge because of the branch cuts (discontinuities) that the zeroes of the kernel might have in the right half-plane of the complex $s_1$-plane. The running example from Section \ref{subs5.1} can thus be considered a simple instance of the inversion algorithm.

\vspace{.2cm}
\paragraph{\textbf{Numerical inversion}}
For the purpose of inverting (\ref{original transform}), consider the Laplace transform of the bivariate tail probability of the waiting time. By a straight-forward integration by parts, this can be related to the Laplace-Stieltjes transform of the waiting time/survival function:

\begin{equation}\label{pony trick}\iint e^{-s_1 u_1-s_2 u_2}\, \prob(W^{(1)}\!>\!u_1,W^{(2)}\!>\!u_2)\,\de u_1\de u_2 =\frac{1}{s_1 s_2}[1- \psi(s_1,0) - \psi(0,s_2)+\psi(s_1,s_2)]. \end{equation}

The key remark is that under mild conditions, this transform is continuous up to the boundary of the non-negative quadrant, as opposed to the Laplace transform of the survival function:

\[\iint e^{-s_1 u_1-s_2 u_2}\, \prob(W^{(1)}\!\leq \!u_1,W^{(2)}\!\leq\!u_2)\,\de u_1\de u_2 \]
 which has by definition a singularity at $(s_1,s_2)=(0,0)$. It is easy to see that for example

\[\lim_{\stackrel{s_1\rightarrow 0}{\mathcal Re\,s_1>0}} \frac{1}{s_1 s_2}[1- \psi(s_1,0) - \psi(0,s_2)+\psi(s_1,s_2)] = \frac{1}{s_2} \left[\frac{\partial \psi}{\partial s_1}(0,s_2) - \frac{\partial \psi}{\partial s_1}(0,0) \right], \]
and even further

\[\lim_{\stackrel{s_1,s_2\rightarrow 0}{\mathcal Re\,s_1,s_2>0}} \frac{1}{s_1 s_2}[1- \psi(s_1,0) - \psi(0,s_2)+\psi(s_1,s_2)] = \frac{\partial^2 \psi}{\partial s_1 \partial s_2}(0,0), \]
and it is clear that this mixed derivative is equal to $\mean [W^{(1)}W^{(2)}]$, which is the same as the left-hand side of (\ref{pony trick}) evaluated at $s_1\!=\!s_2\!=\!0$.
The partial derivatives on the right-hand side above must be considered as limits from the interior of the positive quadrant. The Laplace transform of the ruin function is continuous up to the boundary given that the above partial derivatives exist.

%
%

 The main point of the above discussion is that we may now use the standard form of the multidimensional inversion algorithm developed in ~\cite{denIseger}, for which it is essential that the Laplace transform is regular and continuous up to the boundary of the positive quadrant. The above trick of passing to tail probabilities thus frees one from considering modifications of the inversion algorithm for non-smooth functions (see ~\cite{denIseger}). Once the tail probability/ruin function has been obtained, the survival function follows from an identity similar to (\ref{pony trick}):

 \[\prob(W^{(1)}>x_1, W^{(2)}>x_2) = 1 - \prob(W^{(1)}\leq x_1) -  \prob(W^{(2)}\leq x_2) + \prob(W^{(1)}\leq x_1,\,W^{(2)}\leq x_2),    \]
for any $x_1,\, x_2\geq 0$.

There are no regularity problems with the transforms we will be working with throughout this section because they are all meromorphic in both arguments for positive real values (in some cases they are constructed from branches of various locally meromorphic functions via analytic continuation -- see Example 2).

Above we discussed how to consider the input for the inversion algorithm; some remarks are also needed about the output. This is an $M_1\times M_2$ matrix, that  represents the values of the ruin function $\prob(\tau_\vee(\cdot,\cdot)<\infty)$ on a grid: the entry $(k,l)$ stands for

\[\prob(\tau_\vee((k-1)\Delta_1,(l-1)\Delta_2)<\infty)= \prob(W^{(1)}>(k-1)\Delta_1,W^{(2)}>(l-1)\Delta_2),\]
where $\Delta_i$ are division sizes. The values of the inverted transform are plotted in the figures below for various examples.

\vspace{.2cm}

\noindent\textbf{Example 1.} This is the running example that started at (\ref{running example}).
   We have calculated the one-parameter Wiener-Hopf factor for this example in (\ref{one-parameter WH example}). For the LST of the survival function we need also the positive Wiener-Hopf factor of the projected kernel:

\[ \tilde K(s_1,s_1) = K(s_1,0) = \frac{-9 s_1 - 35 s_1^2 - s_1^3 + 18 s_1^4 +8 s_1^5 + s_1^6}{(-1 + s_1)^2 (2 + s_1)^2 (3 + s_1)^2}.\]
The numerator can be factorized as $s_1(s_1^2+4s_1+1)(s_1^3+4s^2_1+s_1-9)$  where the order 2 polynomial further factorizes as
\[ s_1^2+4s_1+1 =(s_1 +2 -\sqrt{3})(s_1+2+\sqrt{3}) \]
and the zeroes of the factor $s_1^3+4s^2_1+s_1-9$ are

{\small
   \[ v_0 := -\frac{4}{3}  + \frac{1}{3\sqrt[3]{2}}\left(151 - 9 \sqrt{173}\right)^{1/3} + \frac{1}{3\sqrt[3]{2}}\left(151 + 9 \sqrt{173}\right)^{1/3}, \]

   \[v_1:= -\frac{4}{3} - \frac{1}{6\sqrt[3]{2}} (1 + i \sqrt{3}) \left(151 - 9 \sqrt{173}\right)^{1/3} -
 \frac{1}{6\sqrt[3]{2}} (1 - i \sqrt{3}) \left(151 + 9 \sqrt{173}\right)^{1/3}, \]
   \[  v_2 = -\frac{4}{3} - \frac{1}{6\sqrt[3]{2}} (1 - i \sqrt{3}) \left(151 - 9 \sqrt{173}\right)^{1/3} -
 \frac{1}{6\sqrt[3]{2}} (1 + i \sqrt{3}) \left(151 + 9 \sqrt{173}\right)^{1/3},\;\; \mbox{with }v_2=\bar v_1.   \]
}
The positive Wiener-Hopf factor of the projected kernel becomes, cf. (\ref{proj factor}),

\begin{equation}\label{proj factor ex1} K_{pr}^+(s_1) = \frac{(s_1^2+4s_1+1)(s_1-v_1)(s_1-\bar v_1)}{(s_1+2)^2(s_1+3)^2}. \end{equation}
(\ref{one-parameter WH example}) and (\ref{proj factor ex1}) are the necessary components for constructing the LST of the survival function in this example, as given by (\ref{original transform}).
The atom at $(0,0)$ is

\[\prob(W^{(1)}=0,W^{(2)}=0) = \prob(W^{(1)}=0)=K^+_{pr}(0)= \frac{|v_1|^2}{36}\approx 0.204.\]
Here we used the ordering between $W^{(1)}$ and $W^{(2)}$. Below is given the final formula for the transform of the survival function in the original coordinates:

 {\footnotesize
 \[
   \psi(s_1,s_2)= \frac{[(s_1+3)(2s_1+1)+\sqrt{(3s_1+9)(1+3s_1)}][(s_1+3)(2s_1+1)+\sqrt{(s_1+3)(3s_1+13)} ]}{[(s_1+3)(2s_1+2s_2+1)+\sqrt{(3s_1+9)(1+3s_1)}][(s_1+3)(2s_1+2s_2+1)+\sqrt{(s_1+3)(3s_1+13)} ]}\cdot \]

   \vspace{0.2cm}
   \begin{equation}\label{psi rationality}  \frac{|v_1|^2}{36}\frac{(s_1+s_2+2)^2 (s_1+3)^2}{(s_1^2+4s_1+1)(s_1-v_1)(s_1-\bar v_1)}.    \end{equation}
 }

\noindent  $\sqrt{s_1+3}$ cannot be simplified above. The reason is the same as the one given in Remark \ref{branchcuts} below. The marginal transforms are
\[ \psi(s_1,0)=  \frac{|v_1|^2}{36}\frac{(2+s_1)^2(3+s_1)^2}{(s_1^2+4s_1+1)(s_1-v_1)(s_1-\bar v_1)},  \]

\[\psi(0,s_2) = \frac{(3+\sqrt{39})(s_2 +2)^2 }{4( 6s_2+3+\sqrt{39})(s_2+1)}. \]
The atom at zero of $W^{(2)}$ is approximately $0.575$.
We have carried out the numerical inversion for the transform in Example 1: the division size is chosen $\Delta_1=\Delta_2=.1$ and the grid size is $M_1=M_2=2^6$. An important performance measure is the $5\%$ quantile curve of the tail probability -- the two dimensional version of the $5\%$ quantile, also known as value at risk. This is the (not necessarily continuous in general) curve that contains all $(u_1,u_2)$, such that $\prob(W^{(1)}\!>\!u_1,W^{(2)}\!>\!u_2)\geq .05$ and $\prob(W^{(1)}\!>\!u_1+,W^{(2)}\!>\!u_2+)\leq.05.$
To put it simply, the ruin function $\prob(\tau_\vee(\cdot,\cdot)<\infty)$ (see Formula (\ref{naive time to ruin}) and Remark \ref{Remark scaling}) is less than $5\%$ whenever it is evaluated at a point which lies outside the region bounded by this curve in the non-negative quadrant. This, together with several other quantile curves, is displayed in Figure \ref{qcurve ex1} below.

Finally, as a verification, we estimated the ruin function using simulation. Upon choosing suitable bin sizes that account for the atom at $(0,0)$ of $(W^{(1)},W^{(2)})$, the uniform distance between the output of the inversion algorithm and the simulated ruin function is of the order of $10^{-3}$. If we denote with $R(x_1,x_2)= \prob(\tau_\vee (x_1,x_2)<\infty)$, the joint ruin function, the simulation comparison is in Table \ref{tableone}.


\begin{table}[!htbp]
\begin{center}
\begin{tabular}{c ccc ccc cccc}
\hline\hline
 $(x_1,x_2)$  &$(0,0)$ &(2,0) &(2,2) &$(4,0)$ &$(4,2)$ &$(4,4)$ & $(6,0)$ & $(6,2)$ & $(6,4)$ &(6,6) \\ \hline
$\hat R(x_1,x_2)$  &.423  &.297  &.060  &.180  &.050  & .008 & .107 & .034 & .007 & .001\\
$ R(x_1,x_2)$ &.424 &.301 &.060  &.184  &.050 &.008 & .110 & .035 & .007 & .001\\

\hline\hline
\end{tabular}
\end{center}
\caption{Comparison between the simulated ruin function ($\hat R$) and the inverted function ($R$) for various values of the initial capital $(x_1,x_2)$.}
\label{tableone}
\end{table}

\begin{center}
\begin{figure}[htb!]
\includegraphics{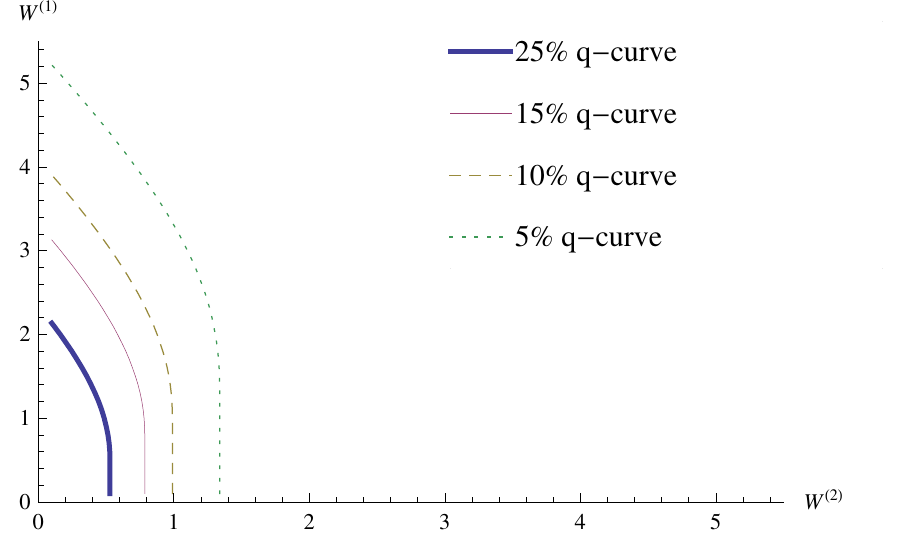}
\caption{ $25\%$, $15\%$, $10\%$, and respectively $5\%$-quantile curves for the ruin function in Example 1. The abscissa corresponds to the values at risk in the second insurance line/the marginal tail of $W^{(2)}$.}
\label{qcurve ex1}
\end{figure}
\end{center}

\begin{Remark}
The quantile curve plots in Figure \ref{qcurve ex1} contain lines which below the main diagonal are straight. This is a consequence of the ordering assumption on the claims (and implicitly on the waiting times) since we can write for $x_1\leq x_2$

\[ \prob(W^{(1)}>x_1,W^{(2)}>x_2) = \prob(W^{(2)}>x_2),\]
because $W^{(2)}>x_2$ implies $W^{(1)}>x_1$ for all $x_1\leq x_2$.
\end{Remark}

\noindent\textbf{Example 2.} The parameters are the same as in Example 1, except now the order is $n=3$. The kernel is

\[\tilde K(s_1,z) = 1- \frac{72}{(3 + s_1)^3 (1 - z)^3 (2 + z)^3} - \frac{12}{(3 + s_1)^2 (1 - z)^2 (2 +
    z)^2} - \frac{2}{(3 + s_1) (1 - z) (2 + z)}.\]
Below we list the zeroes of the numerator. The radicals are again defined when the cut is taken along the negative half of the real axis.

{\small
\[ v_1(s_1) = \frac{-(s_1+3)^2- \sqrt{(s_1+3)^4 -
  4 (s_1+3)^2 (-24 - 14 s_1 - 2 s_1^2 + 2 \sqrt{2} \sqrt{-(3 + s_1)^2})} }{2 (s_1+3)^2}, \]
\[ v_2(s_1) = \frac{-(s_1+3)^2 + \sqrt{(s_1+3)^4 -
  4 (s_1+3)^2 (-24 - 14 s_1 - 2 s_1^2 + 2 \sqrt{2} \sqrt{-(3 + s_1)^2})} }{2 (s_1+3)^2}, \]
\[ v_3(s_1) = \frac{-(s_1+3)^2 - \sqrt{(s_1+3)^4 -
  4 (s_1+3)^2 (-24 - 14 s_1 - 2 s_1^2 - 2 \sqrt{2} \sqrt{-(3 + s_1)^2})} }{2 (s_1+3)^2}, \]
\[ v_4(s_1) = \frac{-(s_1+3)^2 + \sqrt{(s_1+3)^4 -
  4 (s_1+3)^2 (-24 - 14 s_1 - 2 s_1^2 - 2 \sqrt{2} \sqrt{-(3 + s_1)^2})} }{2 (s_1+3)^2}, \]
\[v_5(s_1) =\frac{-(s_1+3) - \sqrt{3} \sqrt{(s_1+3)(1+3s_1)}}{2(s_1+3)},\;\;v_6(s_1) =\frac{-(s_1+3) + \sqrt{3} \sqrt{(s_1+3)(1+3s_1)}}{2(s_1+3)}.\]
}

\begin{Remark}\label{branchcuts} The above formulae cannot be simplified. When choosing a branch for the square root as a function of a complex variable, one has for $a\neq b$:
\[ \sqrt{(z-a)(z-b)}\neq \sqrt{z-a}\sqrt{z-b}. \]
In addition, the term $\sqrt{-(3+s_1)^2}$ is discontinuous (its discontinuity line is $\mathcal Im\,s_1=0$) and it contributes towards the discontinuities of the zeroes $v_i(s_1)$, $i=\overline{1,4}$. 
\end{Remark}

It is not clear a priori which one of the four zeroes to choose when constructing the one-parameter factor $\tilde K_{s_1}(z)$ from (\ref{positive W-H factor}) because, in contrast to Example 1, the branch cuts of $v_1(s_1)$ up to $v_4(s_1)$  cross inside the right half-plane, and  the zeroes $v_i(s_1)$ jump from positive to negative real values when the argument passes between the regions bounded by the cuts in $\mathcal Re\;s_1>0$.

The key observation is that $v_1(s_1)$ is an analytic continuation of $v_2(s_1)$, and $v_3(s_1)$ is an analytic continuation for $v_4(s_1)$.
Therefore, in order to obtain $\tilde K_{s_1}(z)$, one has to glue together (using analytic continuation) the negative branches of $v_1(s_1)$ and $v_2(s_1)$ on the one hand, and of $v_3(s_1)$ and $v_4(s_1)$ on the other, for $\mathcal Re\,s_1\!>\!0$. $v_5(s_1)$ is negative for any $\mathcal Re\,s_1\!>\!0$ so it will always enter the formula for $K^+_{s_1}(\cdot)$ as opposed to $v_6(s_1)$ which is positive and doesn't play any role. Moreover the branch cuts of $v_1(s_1)$ up to $v_4(s_1)$ partition the complex half-plane in 4 regions symmetric around the real axis and the cuts are pairwise parallel lines at angles $\pm\pi/4$. To be more precise, we have to use the 3 different branches of $\tilde K^+_{s_1}(z)$:

\[\tilde K^{1,3}_{s_1}(z)= \frac{(z-v_1(s_1))(z-v_3(s_1))(z-v_5(s_1))}{(z+1)^3}, \,\,\,\, \tilde K^{2,3}_{s_1}(z)= \frac{(z-v_2(s_1))(z-v_3(s_1))(z-v_5(s_1))}{(z+1)^3}, \]

\begin{equation}\label{branches K} \tilde K^{2,4}_{s_1}(v)= \frac{(z-v_2(s_1))(z-v_4(s_1))(z-v_5(s_1))}{(z+1)^3}.\end{equation}
The branches  $C^{1,3}(s_1),$ $C^{2,3}(s_1)$ and $C^{2,4}(s_1)$ are obtained similarly because from (\ref{C(s) formula}), these are related to the corresponding  branches of $\tilde K^+_{s_1}(z)$ by setting $z\!=\!s_1$; the branch cuts and the partition of the complex plane are therefore the same. Since both $C(\cdot)$ and $K^+_{s_1}(\cdot)$ enter Formula (\ref{original transform}), the expression for the LST of the survival function/joint waiting time (see Remark \ref{queueing remark}) is obtained by patching together (via analytic continuation) the positive branches (in the $s_1$-plane) of the generalized Wiener-Hopf factors from Proposition \ref{Proposition W-H 1 parameter} and Theorem \ref{Main result}.

\begin{center}
\begin{figure}[htpb!]
\includegraphics[scale=.7,trim= 30 10 30 20]{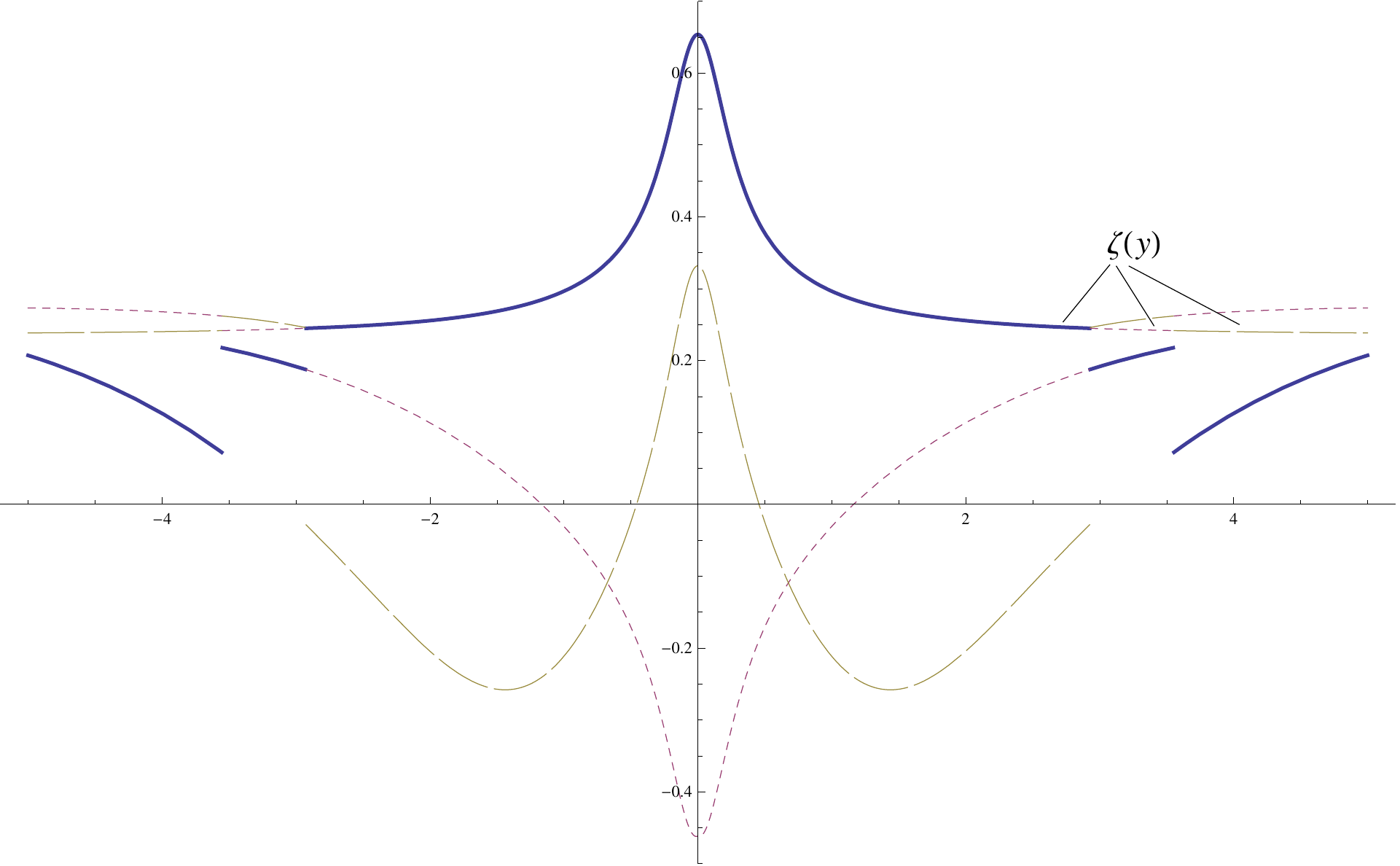}
\caption{The plot of the three branches of the section $\zeta(y)$ from (\ref{example sectioning}).}
\label{branch_example2}
\end{figure}
\end{center}

In Figure \ref{branch_example2} below we plot a section in the three branches of the real part of the LST of the survival measure in Example 2. More precisely consider the section

\begin{equation}\label{example sectioning}\zeta(y):=\mathcal Re\, \psi(iy,14+iy)\end{equation}
that runs along the imaginary axis in the first argument $s_1$ (the argument that generates the discontinuities). From this figure it becomes apparent how the three different branches of $\zeta(y)$ are continuations of each other: the central branch belongs to $\mathcal Re\, \psi^{1,3}(i y, 14+iy)$ (the blue curve). This is continued by the branch $\mathcal Re\,\psi^{2,3}(i y, 14+iy)$ (the dashed red curve) which in turn is continued by $\mathcal Re\, \psi^{2,4}(i y, 14+iy)$ (the orange curve segment) towards the ends of the plot.

\begin{figure}[htb!]
\centering
\subfigure[The joint ruin function /the bivariate tail of the waiting time from Example 2.]{
\centering
\includegraphics[scale=.7]{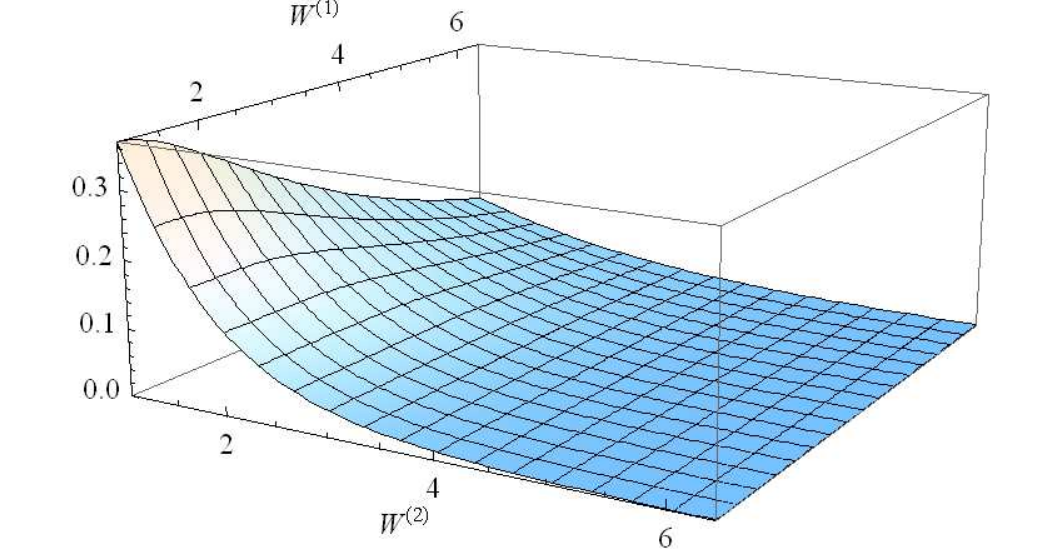}
\label{Fig biv tail}
}
\subfigure[$25\%$, $15\%$, $10\%$, $5\%$, respectively $1\%$-quantile curves for the ruin function in Example 2. The $x$-axis corresponds to the marginal tail of $W^{(2)}$.]{
\centering
\includegraphics[scale=.7]{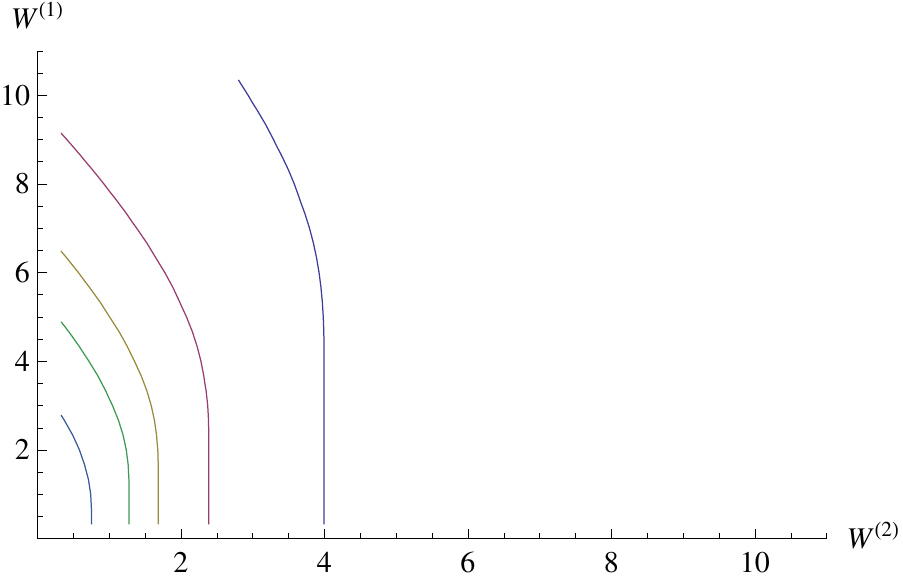}
\label{Fig quantile ex2}
}
\caption{}
\end{figure}

We numerically inverted the expression obtained for the transform $\psi(s_1,s_2)$ using den Iseger's algorithm~\cite{denIseger}. Again, the division size is $\Delta = .1$ and the grid size is $M = 2^6$ in both directions.
In Figure \ref{Fig biv tail} the plot of the ruin function $\prob(\tau_\vee(\cdot,\cdot)<\infty)$ is presented (see Formula (\ref{naive time to ruin}) and Remark \ref{Remark scaling}) or equivalently the tail of the equilibrium distribution for the bivariate waiting time $(W^{(1)},W^{(2)})$ (Remark \ref{queueing remark}) in Example 2.
We can write in general

\[\prob(W^{(1)}>0,W^{(2)}>0)= 1  - \prob(W^{(1)}=0) - \prob(W^{(2)}=0) + \prob(W^{(1)}=0,W^{(2)}=0), \]
and because of the ordering, $\prob(W^{(1)}=0,W^{(2)}=0)= \prob(W^{(1)}=0).$ And then the value at $(0,0)$ of the joint ruin function is $\prob(W^{(2)}>0)$, which in this example approximately equates $0.37$.

Finally, in Figure \ref{Fig quantile ex2}  we present various quantile curves for the ruin function/stationary tail of the waiting time.

\begin{figure}
\centering
\subfigure[25$\%$ quantile curves]{
\includegraphics[scale=.7]{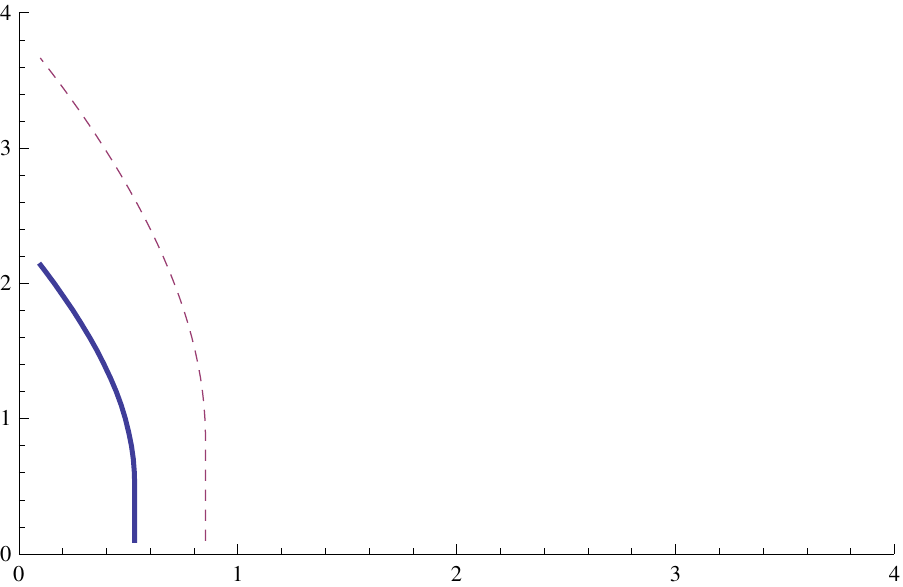}
}
\subfigure[15$\%$ quantile curves]{
\includegraphics[scale=.7]{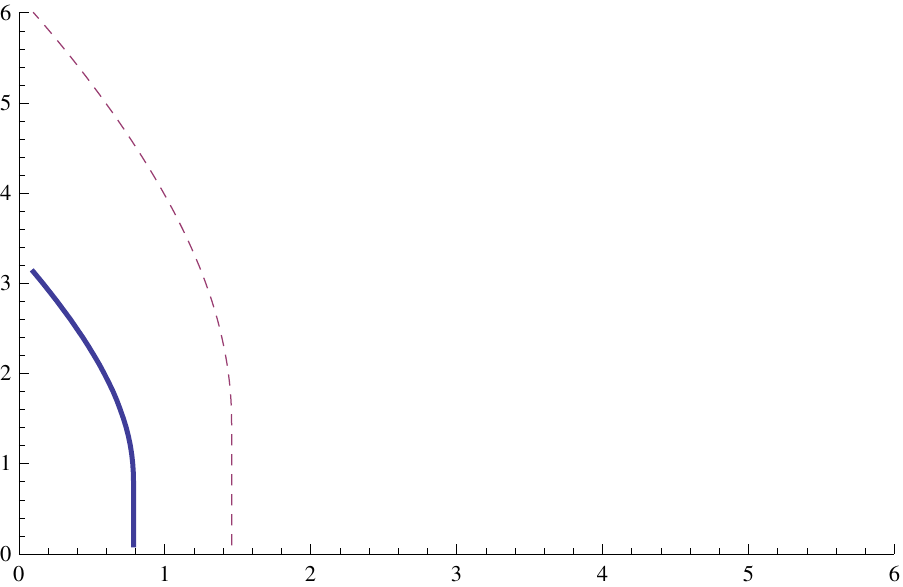}
}
\subfigure[10$\%$ quantile curves]{
\includegraphics[scale=.7]{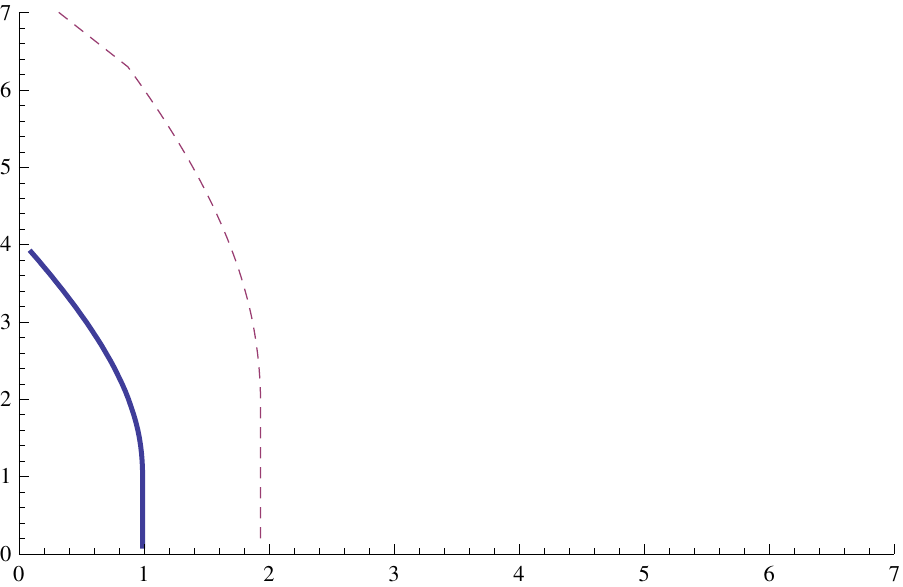}
}
\subfigure[5$\%$ quantile curves]{
\includegraphics[scale=.7]{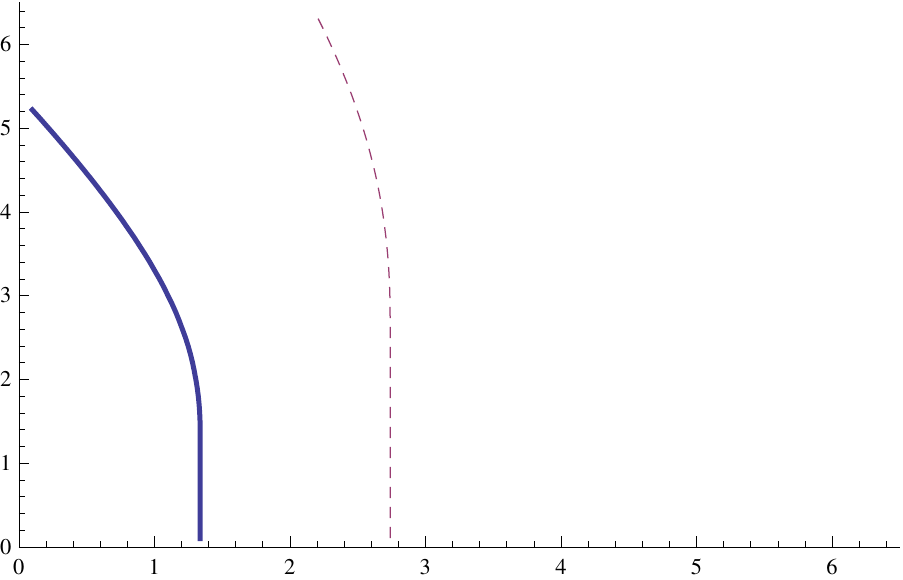}
}
\caption{Comparison of risks. Dashed curves correspond to decoupled input.} \label{figure comparison}
\end{figure}

\vspace{.2cm}

\textbf{Comparison of risks}. In Figure \ref{figure comparison}  we compare the results for quantile curves of the ruin function of Example 1 with  the quantile curves of the ruin function for the case where the input is decoupled. By this we mean we take three samples $N_1$, $N_2$ and $N_3$ from a uniform distribution on $\{1,2\}$ and define the random variables $A_{dec}= \sum_{i=1}^{N_1} A_i$, $D_{dec}= \sum_{i=1}^{N_2} D_i$, $B^{(2)}_{dec}= \sum_{i=1}^{N_3} B^{(2)}_i$, where $(A_i)_{i\leq n}$, $(B^{(2)}_i)_{i\leq n}$, and $(D_i)_{i\leq n}$ are mutually independent sequences of exponential random variables with rates $\lambda$, $\mu$ and $\mu_D$ respectively. In this case the inter-arrival time becomes independent of the claim size vector, while marginally $A_{dec}$, $B^{(2)}_{dec}$ and $B^{(1)}_{dec}$ have the same distribution as in Example 1. The kernel for this instance is

\[  \tilde K(s_1,z) = \left(\frac{9}{2 (3 + s_1)^2} + \frac{3}{2 (3 + s_1)}\right) \left(\frac{1}{2 (1 - z)^2} + \frac{1}{
   2 (1 - z)}\right) \left(\frac{2}{(2 + z)^2} + \frac{1}{2 + z}\right).\]

The zeroes of the numerator as a polynomial in $z$ are already too complicated to present here. This instance is similar to Example 2 in terms of the analytic behaviour of these zeroes.

The main point is that, similarly as in \cite{paper1}, numerical data  suggests that the ruin functions corresponding to positively correlated input on the one hand, and the ruin functions for decoupled input on the other are stochastically ordered (Figure \ref{figure comparison}).

\vspace{.3cm}
\noindent\textbf{Example 3 (proportional reinsurance)}
This is the case with proportional claims. There is a common arrival process such that the inter-arrival time $A_n$ is correlated with the claim size $B_n$, and $\alpha B_n$ is deducted from the first insurance line and $(1-\alpha)B_n$ from the second.

We take $N\sim$ Unif$\{1,2,3\}$, $\lambda=\mu=1$, $\alpha=3/4$ and unit income rates.
For the purpose of comparing risks, we will consider three instances for the random vector $(A,B)$:

\begin{align*}
 &\mbox{positive correlation: } (A,B)_{pos}\sim(\mbox{ Erlang}(N,\lambda),\mbox{ Erlang}(N,\mu)), \\
 &\mbox{independence: }(A,B)_0\sim(\mbox{Erlang}(N_1,\lambda),\mbox{ Erlang}(N_2,\mu)),\mbox{ with }N_1,\; N_2 \mbox{ two copies of }N. \\
 &\mbox{negative correlation: }(A,B)_{neg}\sim\,(\mbox{Erlang}(N,\lambda),\mbox{ Erlang}(4-N,\mu)).
\end{align*}
The kernels corresponding to these instances are

\[ \tilde K_{pos}(s_1,z) = \frac{1}{3(1-z)(1+s_1/2 + z/4)} + \frac{1}{3(1-z)^2(1+s_1/2 + z/4)^2} + \frac{1}{3(1-z)^3(1+s_1/2 + z/4)^3}, \]

\begin{align*} \tilde K_{0}(s_1,z) = &\left(\frac{1}{3(1-z)} + \frac{1}{3(1-z)^2} + \frac{1}{3(1-z)^3}\right)\left(\frac{1}{3(1+s_1/2 + z/4)} + \frac{1}{3(1+s_1/2 + z/4)^2} \right. \\ &+\left.\frac{1}{3(1+s_1/2 + z/4)^3}\right),\end{align*}

\[ \tilde K_{neg}(s_1,z) = \frac{1}{3(1-z)(1+s_1/2 + z/4)^3} + \frac{1}{3(1-z)^2(1+s_1/2 + z/4)^2} + \frac{1}{3(1-z)^3(1+s_1/2 + z/4)}. \]

These functions stand for $\mean e^{-[s_1(2\alpha-1) + z (1-\alpha)]B +zA }$ under the three couplings. The correlations between the variables $A$ and $B$ can also be read directly from the shapes of these transforms.

%
%
%

\begin{figure}[htb!]
\centering
\subfigure[]{
\includegraphics[scale=.41]{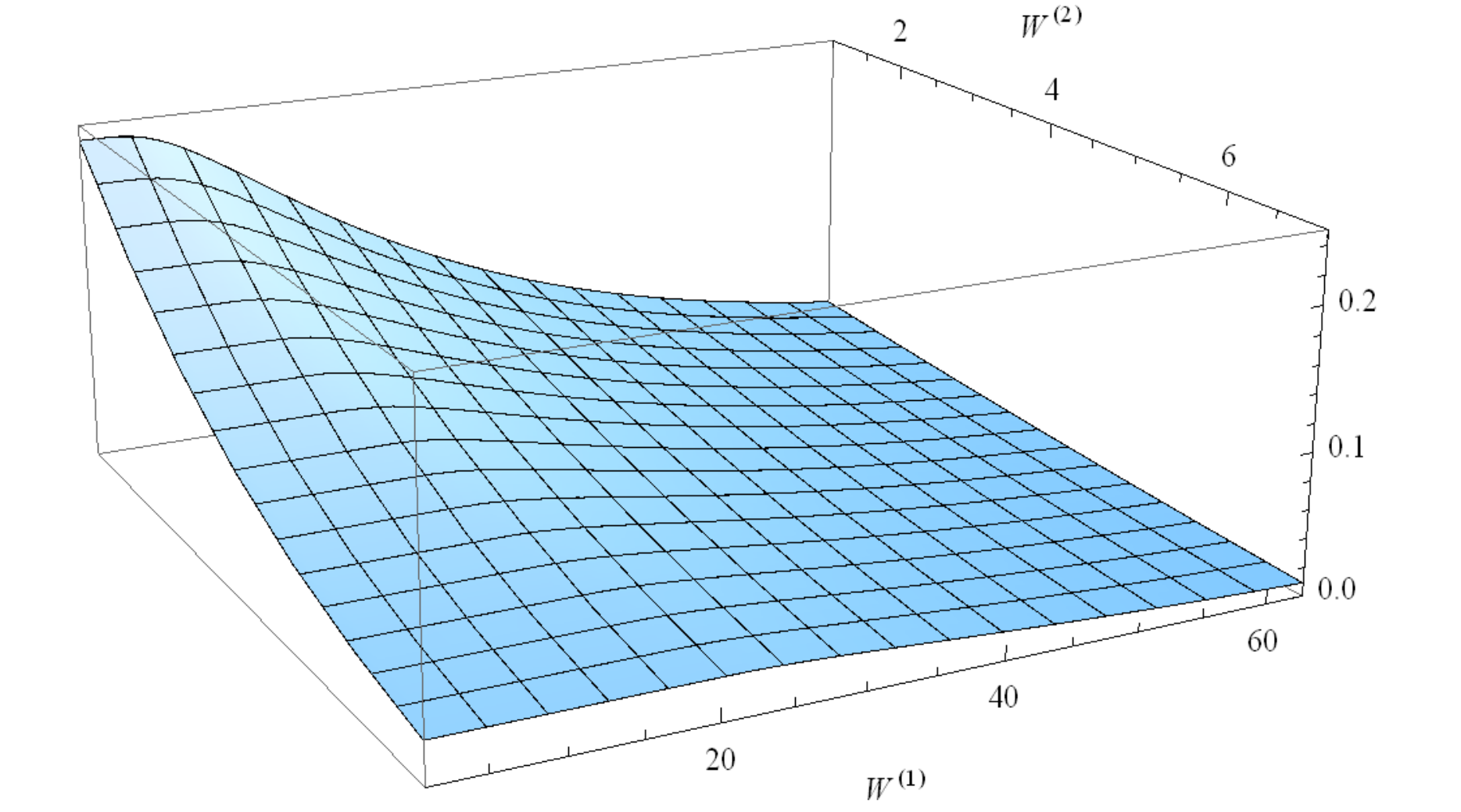}
}
\subfigure[]{
\includegraphics[scale=.7]{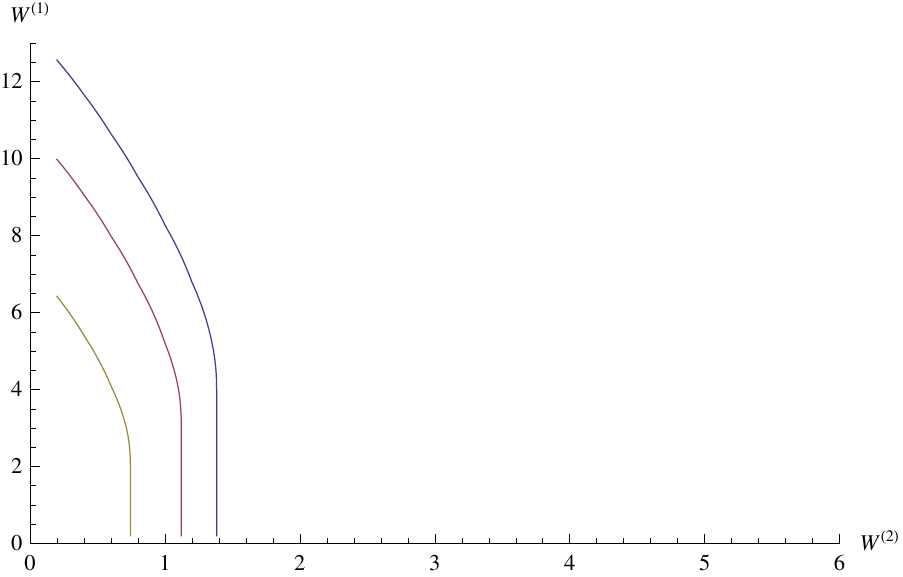}\label{qcurves_prop}
}
\caption{Bivariate tail (left) and respectively 10$\%$, 5$\%$, and 3$\%$ quantile curves (right) for proportional reinsurance with negative correlation.} \label{Fig prop}
\end{figure}

 Numerical illustrations are in Figure \ref{Fig prop}.

\begin{figure}
\centering
\subfigure[3$\%$ quantile curves]{
\includegraphics[scale=.7]{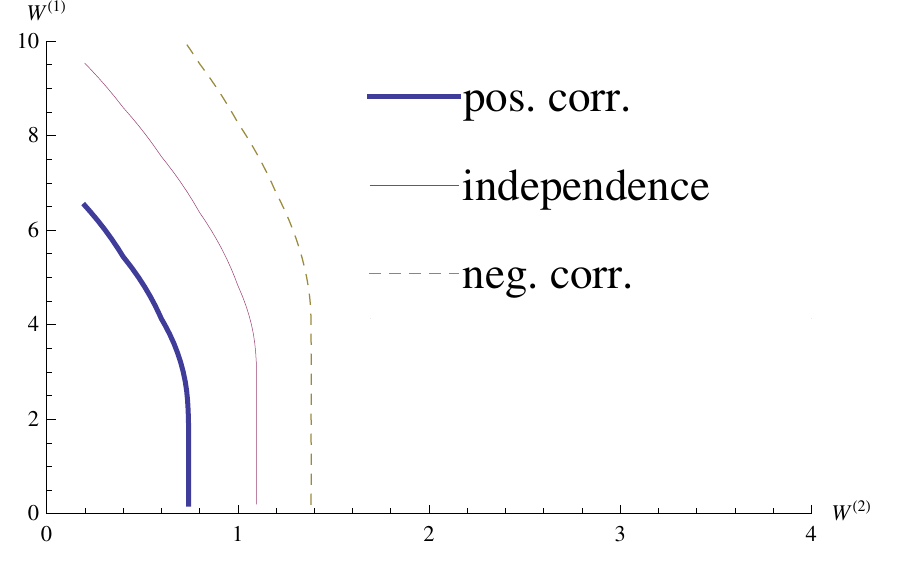}
}
\subfigure[5$\%$ quantile curves]{
\includegraphics[scale=.7]{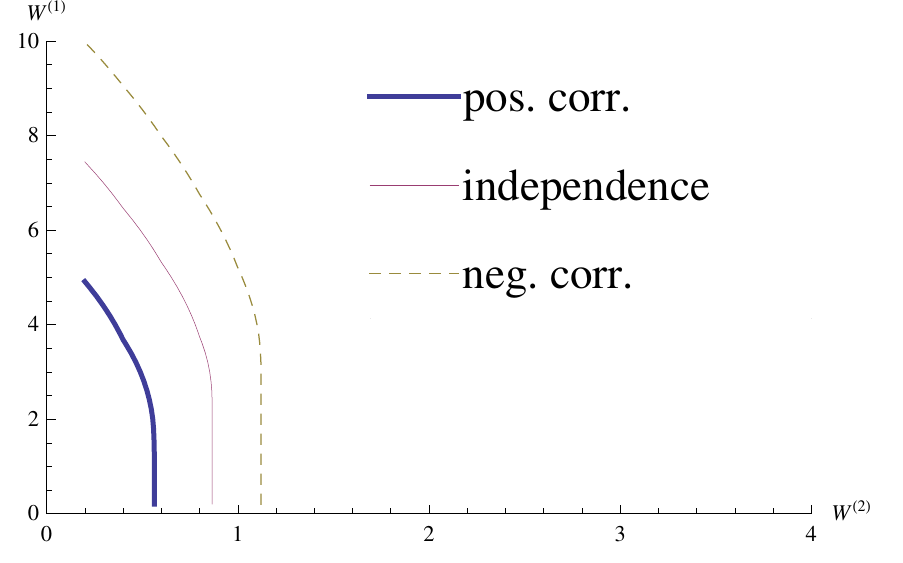}
}
\subfigure[10$\%$ quantile curves]{
\includegraphics[scale=.7]{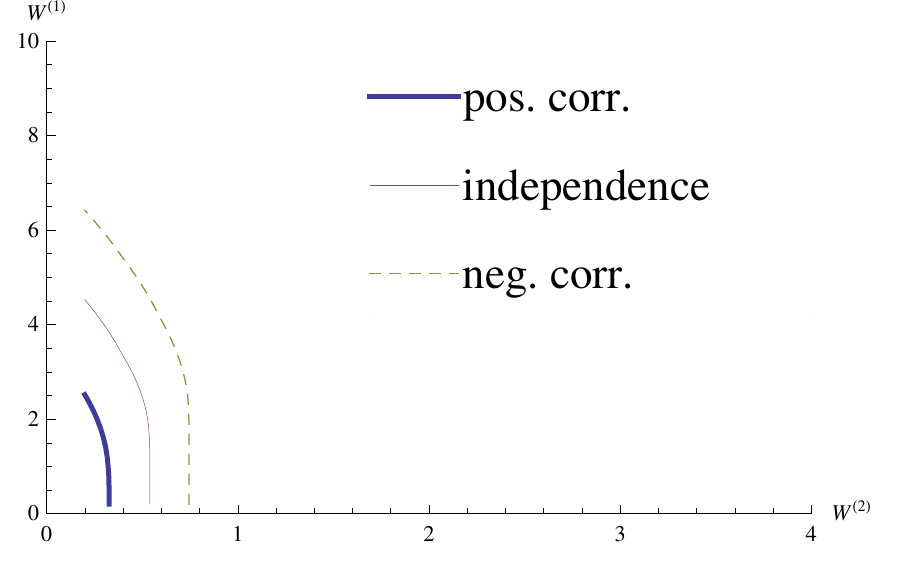}
}
\caption{Comparison of risks.} \label{figure prop comparison}
\end{figure}

\textbf{Comparison of risks}
In the table below we present various points at which a specific ruin probability is achieved. Given a fixed value at risk, any value for starting
capital $(x_1,x_2)$ lying on the respective quantile curve will achieve this. Interestingly, the risks are ordered between the various types of correlations (see also Figure \ref{figure prop comparison}). 
Let us denote by $R(x_1,x_2)=\prob(\tau_\vee(x_1,x_2)<\infty),$ the probability that eventually, both lines are ruined. Then for the three types of correlation we give the values of the function $R$ in Table \ref{tabletwo}. Positive correlation gives the lowest ruin probabilities for any starting capital considered.

\begin{table}[!htbp]
\begin{center}
\begin{tabular}{c ccc ccc cccc}
\hline\hline
 $(x_1,x_2)$  &$(0,0)$ &(2.4,0) &(4.8,0) &$(4.8,.4)$ &$(6.4,.4)$ &$(6.4,.8)$ &(9,0.4) & $(9,.8)$ & $(11.8,.8)$  \\ \hline
$ R_{neg}(x_1,x_2)$  &.2388  &.1995  &.1309  &.0862  &.0648  & .0402 & .0397 & .0253 & .0149 \\
$ R_{0}(x_1,x_2)$ & .1922 &.1516 &.0896 &.0536  &.0375  &.0214 &.0203 & .0120 & .0061  \\
$ R_{pos}(x_1,x_2)$ &.1381 &.0979 &.0486  &.0237  &.0148 &.0070 & .0065 & .0033 & .0013 \\
\hline\hline
\end{tabular}
\end{center}
\caption{Comparison between the joint ruin functions $R_{neg}$, $R_0$ and  $R_{pos}$ respectively, for the various types of correlation.}
\label{tabletwo}
\end{table}

\section*{Conclusions}
The working assumptions of rationality and ordering for the trivariate input made up from the generic claim vector together with the preceding inter-arrival time, allow one to obtain detailed numerical results, for the joint ruin probability as a function of the initial reserves. Our numerical results suggest that when comparing ruin functions that correspond to various correlation structures between claim intervals and claim sizes, positive correlation among these cause lower values of the ruin probability compared to zero correlation, and even more so, when compared to negative correlation.

By using the relations described in Section 2, one can recover the other types of ruin/survival functions, where the numerical inversion of the marginal transforms, when needed, can be carried out using the one-dimensional inversion algorithm.

\appendix
\section{}\label{Appendix}

\begin{AProposition}[On the zeroes of the kernel $1-\tilde K(s_1,z)$]\label{Rouche problem} For each $s_1$, $\mathcal Re\;s_1\geq 0$, $g(s_1,z)- f(s_1,z)$ and $g(s_1,z)$ have the same number of zeroes in $\mathcal Re\; z\geq 0.$
\end{AProposition}

\proof We show that $g(s_1,z)$ dominates $f(s_1,z)$ on a suitably chosen contour in the complex $z$-plane. From this the claim in the proposition will follow via Rouch\'e's theorem~\cite{Titchmarsh}, p.116. Consider the contour which is made up from the extended arc:

\[\mathcal C_\epsilon:=\left\{re^{i\varphi};\;\varphi\in \left[-\pi/2-\arccos \epsilon,\pi/2+\arccos\epsilon\right] \right\},\] together with the line segment

\[I:=\left\{-\epsilon+i\omega;\;|\omega|\in\left[0,r\sqrt{1-\epsilon^2} \right]\right\}.\]
The rationality of the transform ensures that $\tilde K(s_1,z)$ can be analytically continued on a thin strip: $\mathcal Re\;z<0$, $|\mathcal Re\;z|<\epsilon$.
  We first consider the contour $\mathcal C_\epsilon$. On the one hand we have the triangle inequality for $f$: $|f(s_1,z)|\leq \bar f(|s_1|,|z|)$, where $\bar f$ is a polynomial with the same degree as $f$. On the other hand, if we represent $g_{s_1}(z)$ for fixed $s_1$ as $g_{s_1}(z) = a_m(s_1)\prod_{i} (z - \xi_i(s_1))$, with $\xi_i(s_1)$ its zeroes and $a_m(s_1)$ the coefficient of $z^m$, $m=m(s_1)=\deg g_{s_1}(z)$,  then the same triangle inequality gives a lower bound for $|g(s_1,z)|$:

  \[|g(s_1,z)| =  |a_m(s_1)| \prod_{i} \left|z - \xi_i(s_1)\right|\geq |a_m(s_1)|\prod_{i} \left(|z| - |\xi_i(s_1)|\right)=: \bar g(s_1,|z|), \] with the remark that $\bar g_{s_1}(|z|)$ has the same degree as $g_{s_1}(z)$.

  Now we can bound for $r$ sufficiently large, such that the interior of $\mathcal C_\epsilon \cup I$ contains all the zeroes of $\bar g_{s_1}(|z|)$ and  $z\in \mathcal C_\epsilon$:

\[ \left|\frac{f(s_1,z)}{g(s_1,z)}\right| \leq  \frac{\bar f(|s_1|,|z|)}{\bar g(s_1,|z|)} \rightarrow 0, \;\;\;\mbox{ as $r\rightarrow \infty$}.  \]

Convergence holds because the degree of the numerator is strictly less than that of the denominator, by Assumption \ref{assumption rational}. This establishes the bound $|g(s_1,z)|> |f(s_1,z)|$ on $\mathcal C_\epsilon$, for $r$ large enough.

For the segment $I$, we use the safety loading condition for the second($!$) line: $c^{(2)} \mean A- \mean B^{(2)} >0.$ That is, we start with the fact $\frac{d}{dz}\frac{f(0,z)}{g(0,z)}|_{z=0}= c^{(2)}\mean A -\mean B^{(2)}>0.$ So for $\epsilon>0$ sufficiently small, $\frac{f(0,-\epsilon)}{g(0,-\epsilon)}<\frac{f(0,0)}{g(0,0)}=1$. Then we can write for $z\in I$:

  \[\left|\frac{f(s_1,\epsilon+i\omega)}{g(s_1,\epsilon+i\omega)}\right| \leq \mean \left(|e^{-s_1 D}|\cdot |e^{-\epsilon(A-B^{(2)})}|\cdot |e^{-i\omega(A- B^{(2)})}|\right) \leq \mathbb Ee^{-\epsilon(A-B^{(2)})}=\frac{f(0,-\epsilon)}{g(0,-\epsilon)}<1. \]
Above we used the rough bound $|e^{-s_1D}|\leq 1$. This completes the proof.

Notice that the key role in the proof is played by $\rho_2<1$ and not $\rho_1<1$ ($\rho_2<\rho_1$).


%
%
%



\begin{thebibliography}{10}

\bibitem{AsmussenRuin}
S.~Asmussen and H.~Albrecher.
\newblock {\em Ruin Probabilities}.
\newblock World Scientific Publ. Cy., Singapore, 2010.

\bibitem{APP1}
 F.~ Avram, Z.~Palmowski, and M.~Pistorius.
\newblock A two-dimensional ruin problem on the positive quadrant.
\newblock {\em Insurance: Mathematics and Economics}, {\bf 42}(1), 227--234, 2008.

\bibitem{APP2}
 F. Avram, Z. Palmowski and  M. Pistorius.
\newblock Exit problem of a two-dimensional risk process from the quadrant:
  {E}xact and asymptotic results.
\newblock {\em Annals of Applied Probability\/} {\bf 18,} 2421--2449, 2008.

\bibitem{Badescu}
A.L.~Badescu, E.C.K.~Cheung, and L. Rabehasaina.
\newblock A two dimensional risk model with proportional reinsurance.
\newblock {\em J. Appl. Prob.}, {\bf 48}, 749--765, 2011.

\bibitem{paper1}
E.S.~Badila, O.J.~Boxma, and J.A.C.~Resing.
\newblock Queues and risk processes with dependencies.
\newblock {\em Stochastic Models}, {\bf 30}(3), 390--419, 2013.

\bibitem{paper2}
E.S.~Badila, O.J.~Boxma, and J.A.C.~Resing.
\newblock Queues and risk models with simultaneous arrivals.
\newblock {\em Adv. Appl. Prob.}, {\bf 46}(3), 812--831, 2014.

\bibitem{Bladt:MVME}
M.~Bladt and B.F. Nielsen.
\newblock Multivariate matrix-exponential distributions.
\newblock {\em Stochastic Models}, {\bf 26}(1), 1--26, 2010.

\bibitem{Cai&Li}
J.~Cai, H.~Li.
\newblock Multivariate risk model of phase type.
\newblock {\em Insurance:Mathematics and Economics}, {\bf 36}, 137--152, 2005.

\bibitem{Chan}
W.-S.~Chan, H.~Yang, and  L.~Zhang.
\newblock Some results on ruin probabilities in a two-dimensional risk model.
\newblock {\em Insurance: Mathematics and Economics,} {\bf 32,} 345--358, 2003.

\bibitem{CohenSingleServer}
J.W. Cohen.
\newblock {\em The Single Server Queue}.
\newblock North Holland, 1982.

\bibitem{BVP}
J.~W. Cohen and O.~J. Boxma,
\newblock {\em Boundary Value Problems in Queueing System Analysis}.
\newblock North-Holland Publ. Cy., Amsterdam, 1983.

\bibitem{denIseger}
P. den Iseger.
\newblock Numerical transform inversion using Gaussian quadrature.
\newblock {\em Probab. Eng. Inform. Sc.}, {\bf 20}, 1--44, 2006.

\bibitem{Gong}
L.~Gong, A.L.~Badescu, and  E.C.K.~Cheung.
\newblock Recursive methods for a multi-dimensional risk process with common shocks.
\newblock {\em Insurance: Mathematics and Economics,} {\bf 50}, 109--120, 2012.

\bibitem{KingmanAlgebra}
J.F.C. ~Kingman.
\newblock On the algebra of queues.
\newblock {\em J. Appl. Prob.,} {\bf 3}(2), 285--386, 1966.

\bibitem{Kulkarni:MPH*}
V.G. Kulkarni.
\newblock A new class of multivariate phase type distributions.
\newblock {\em Operations Research}, {\bf 37}(1), 151--158, 1989.

\bibitem{Prabhu80}
N.U.~Prabhu.
\newblock {\em Stochastic {S}torage {P}rocesses. {Q}ueues, {I}nsurance {R}isk, and {D}ams}.
\newblock Springer Verlag, 1980.


\bibitem{Stoyan&Daley}
D.~Stoyan.
\newblock {\em Comparison Methods for Queues and Other Stochastic Models}.
\newblock John Wiley \& Sons Ltd., 1983.

\bibitem{Sundt}
 B.~Sundt.
 \newblock On multivariate {P}anjer recursions.
\newblock {\em ASTIN Bull.}, {\bf 29}, 29--45, 1999.

\bibitem{Titchmarsh}
E.C. Titchmarsh.
\newblock {\em The Theory of Functions}.
\newblock Oxford University Press, 2nd edition, 1939.




\end{thebibliography}
\end{document}